\documentclass{amsart}

\input epsf

\usepackage{amsfonts,amsthm,amsmath,amssymb,latexsym}
\usepackage{graphicx,color}
\usepackage[all]{xy}

\def\D{{\bf D}}

\def\S{\mathcal{S}^p}

\begin{document}

\author{Vladimir Markovic and Dragomir \v Sari\' c}

\address{Department of Mathematics, Stony Brook University,
Stony Brook, NY 11794-3651} \email{markovic@math.sunysb.edu}

\address{Department of Mathematics, Queens College of CUNY,
65-30 Kissena Blvd., Flushing, NY 11367}
\email{Dragomir.Saric@qc.cuny.edu}

\title[Small elements in the Modular group]
{The Teichm\"uller distance between finite index subgroups of
$PSL_2(\mathbb{Z})$}


\date{\today}

\maketitle

\begin{abstract}
For a given $\epsilon >0$, we show that there exist two finite index
subgroups of $PSL_2(\mathbb{Z})$ which are
$(1+\epsilon)$-quasisymmetrically conjugated and the conjugation
homeomorphism is not conformal. This implies that for any
$\epsilon>0$ there are two finite regular covers of the Modular once
punctured torus $T_0$ (or just the Modular torus) and a
$(1+\epsilon)$-quasiconformal between them that is not homotopic to
a conformal map. As an application of the above results, we show
that the orbit of the basepoint in the Teichm\"uller space $T(\S )$
of the punctured solenoid $\S$ under the action of the corresponding
Modular group (which is the mapping class group of $\S$ \cite{NS},
\cite{Odd}) has the closure in $T(\S )$ strictly larger than the
orbit and that the closure is necessarily uncountable.
\end{abstract}

\section{Introduction}

Let $F$ be a quasiconformal map between two Riemann surfaces. By

$$
\mu(F)={{\bar {\partial} F } \over {\partial{F}}} ,
$$

we denote the Beltrami dilatation (or just the dilatation) of  $F$.
The function

$$
K(F)={{1+|\mu|}\over{1-|\mu|}},
$$
\noindent is called the distortion function of $F$. If $K \geq 1$ is
such that $1 \leq K(F) \leq K$ a.e. we say that $F$ is
$K$-quasiconformal. If $F$ is  homeomorphism of the unit disc onto
itself  that is $1$-quasiconformal then $F$ is a M\"obius
transformations. \vskip .3 cm

Let $f:{\bf S}^{1} \to {\bf S}^{1}$ be a homeomorphism of the unit
circle onto itself. We say that $f$ is a $K$-quasisymmetric map if
there exists a $K$-quasiconformal map $F:{\bf D} \to {\bf D}$ of the
unit disc onto itself so that the continuous extension of $F$ on
${\bf S}^1$ agrees with $f$ (recall that every quasiconformal maps
of the unit disc onto itself extends continuously to a homeomorphism
of the unit circle). If $f$ is homeomorphism of the unit circle onto
itself that is $1$-quasisymmetric then $f$ is conformal, that is $f$
is a M\"obius transformations.

\vskip .3cm

The first main result in this paper is

\vskip .3 cm

\paragraph{\bf Theorem 1} {\it For every $\epsilon >0$ there exist two
finite index subgroups of $PSL_2(\mathbb{Z})$ which are conjugated
by a $(1+\epsilon )$-quasisymmetric homeomorphism of the unit circle
and this conjugation homeomorphism is not conformal.}

\vskip .3 cm

Unless stated otherwise by a  Riemann surface we always mean a
Riemann surface of finite type. Every such Riemann surface is
obtained by deleting at most finitely many points from a closed
Riemann surface. If $M$ and $N$ are Riemann surfaces we say that a
map $\pi:N \to M$ is a finite degree, regular covering if $\pi$ is
holomorphic, of finite degree and locally univalent (we also say
that $N$ is a cover of $M$). (Some people prefer the term unbranched
covering instead of regular covering.) Unless stated otherwise all
coverings are assumed to be regular and of finite degree. Note that
if $M$ and $N$ have punctures then the regularity assumption does
not imply that $\pi$ is locally univalent in a neighborhood of a
puncture (which is only natural since the punctures are not part of
the corresponding surface).

\vskip .3 cm

Given two Riemann surfaces $M$ and $N$ the Ehrenpreiss conjecture
asks if for every $\epsilon >0$ there are coverings $M_{\epsilon}\to
M$ and $N_{\epsilon}\to N$ such that there exists a $(1+\epsilon
)$-quasiconformal map $F:M_{\epsilon}\to N_{\epsilon}$. In this case
we say that $M_{\epsilon}$ and $N_{\epsilon}$ are $\epsilon$-close.
It is easy to see that this conjecture is true if $M$ and $N$ are
tori.

\vskip .3cm

Recall that the notion of quasiconformal
$\epsilon$-closeness  between  hyperbolic Riemann surfaces is in
fact a geometric property . After endowing the Riemann surfaces $M$
and $N$ with the corresponding hyperbolic metrics it is well known
that a  $(1+\epsilon )$-quasiconformal map $F:M_{\epsilon}\to
N_{\epsilon}$ is isotopic to a $(1+\delta)$-biLipschitz
homeomorphism such that $\delta\to 0$ when $\epsilon\to 0$ (the
function $\delta=\delta(\epsilon)$ does not depend on the choice of
surfaces $M$ and $N$).  One such biLipschitz map is obtained by
taking the barycentric extension \cite{DE} of the boundary map of
the lift to the universal covering of $f:M_{\epsilon}\to
N_{\epsilon}$ (this observation was made in \cite{Ge}). \vskip .3 cm
There are no known examples of hyperbolic Riemann surfaces $M$ and
$N$, such that for every $\epsilon>0$ there are coverings
$M_{\epsilon}\to M$ and $N_{\epsilon}\to N$ which are
$\epsilon$-close, unless $M$ and $N$ are commensurate. We say that
$M$ and $N$ are commensurate if they  have a common cover (recall
that we assume throughout the paper that all coverings are regular).
If   $M$ and $N$ are commensurate one can say that $M$ and $N$ have
coverings that are $0$-close.

\vskip .3 cm

It is not difficult to see (see the last section for the proof) that
if the Ehrenpreiss conjecture had a positive answer then for any
Riemann surface $M$ and for any $\epsilon >0$, there would exist two
coverings $M_1,M_2\to M$ and a $(1+\epsilon)$-quasiconformal
homeomorphism $F:M_1 \to M_2$ that is not homotopic to a conformal
map. In general, for a given Riemann surface $M$ the problem of
constructing two such covers $M_1$ and $M_2$ and the corresponding
$(1+\epsilon)$-quasiconformal map $F$ (where $F$ is not homotopic to
a conformal map) seems to have a similar degree of  difficulty as
the Ehrenpreiss conjecture.

\vskip .3 cm

\paragraph{\bf Problem} {\it Let $M$ be a hyperbolic Riemann surface. Is
it true that for every $\epsilon >0$ there exist two coverings
$M_1,M_2\to M$ and a  $(1+\epsilon)$-quasiconformal homeomorphism
$F:M_1 \to M_2$ that is not homotopic to a conformal map?}

\vskip .3cm

If such covering surfaces $M_1$ and $M_2$ exist they may be
conformally equivalent. In that case the homeomorphism $F$ is not
allowed to be homotopic to any such conformal equivalence. If $M_1$
and $M_2$ are not conformally equivalent then we could say that the
coverings $M_1$ and $M_2$ are $\epsilon$-close but not $0$-close.

\vskip .3cm In this paper we show that this problem has a positive
answer  for the Modular  torus (and any other Riemann surface
commensurate with it).

\vskip .3cm

\paragraph{\bf Corollary 1}  {\it Let $T_0$ denote the Modular torus. Then
for every $\epsilon >0$ there are finite degree, regular coverings
$\pi_1:M_1 \to T_0$ and $\pi_2:M_2\to T_0$, and a
$(1+\epsilon)$-quasiconformal homeomorphism $F:M_1 \to M_2$ that is
not homotopic to a conformal map.}

\vskip .3cm
\begin{proof} This follows directly from Theorem 1. Assume that $G_1,G_2<
PSL_2(\mathbb{Z})$ are two finite index subgroups that are
conjugated by $(1+\epsilon)$-quasisymmetric map of the unit circle
that is not conformal. Let $M_i$, $i=1,2$, be the Riemann surface
that is conformally equivalent to the quotient $\D/G_i$. Then $M_1$
and $M_2$ satisfy the assumptions in the statement of this
corollary.
\end{proof}

\vskip .3cm

Consider the coverings $\pi_1:M_1 \to T_0$ and $\pi_2:M_2\to T_0$,
where the surfaces $M_1$ and $M_2$ are from  Corollary 1. Then for
$\epsilon$ small enough the $(1+\epsilon)$-quasiconformal map $F:M_1
\to M_2$ can not be a lift of a self homeomorphism of $T_0$. That
is, there is no homeomorphism $\hat{F}:T_0 \to T_0$ so that $\hat{F}
\circ \pi_1=\pi_2 \circ F$. The non-existence of such a map
$\hat{F}$ follows from the discreteness of the action of the Modular
group on the Teichm\"uller space of the surface $T_0$. This
illustrates  what is difficult about proving Corollary 1. An
important ingredient in the proof of Theorem 1. is the fact that
$PSL_2(\mathbb{Z})$ is an arithmetic lattice. At the moment we can
not prove this result for other punctured surfaces (or for any
closed surface, not even those whose covering groups are
arithmetic). However, already from Corollary 1. one can make strong
conclusions about the Teichm\"uller space of the punctured solenoid.

\vskip .3cm

\paragraph{\bf Remark}  {In  \cite{LR} Long and Reid  defined the notion of
pseudo-modular surfaces and have shown their existence. As an important special case,
it would be interesting to examine whether one can prove the above corollary for
a pseudo-modular surface instead of the Modular torus.}

\vskip .3 cm

Recall that the inverse limit $\mathcal{S}$ of the family of all
pointed regular finite covers of a closed hyperbolic Riemann surface
is called {\it the Universal hyperbolic solenoid} (see \cite{Sul}).
It is well known that the commensurator group of the fundamental
group of a closed Riemann surface acts naturally on the
Teichm\"uller space $T(\mathcal{S})$ of the solenoid $\mathcal{S}$
\cite{NS}. Sullivan has observed that the Ehrenpreiss conjecture is
equivalent to the question whether the orbits of this action are
dense in $T(\mathcal{S})$ with respect to the corresponding
Teichm\"uller metric.

\vskip .3 cm

In \cite{Sul} mainly closed Riemann surface have been considered (as
a model how such holomorphic inverse limits should be constructed).
We consider the family of all pointed coverings of some fixed once
punctured torus $T$. The punctured solenoid $\S$ is the inverse
limit of the above family \cite{PS}. The covers of the punctured
torus $T$ are regular, but as we already pointed out, the coverings
can be naturally extended to the punctures in the boundary and  are
allowed to be branched over those boundary punctures. The punctured
solenoid $\S$ is an analog, in the presence of the punctures, of the
universal hyperbolic solenoid $\mathcal{S}$. The (peripheral
preserving) commensurator group $Comm_{per}(\pi_1(T))$ of the
fundamental group $\pi_1(T)$ of $T$ acts naturally on the
Teichm\"uller space $T(\S )$ of the punctured solenoid $\S$. We
consider the orbit space $T(\S )/Comm_{per}(\pi_1(T))$.

\vskip .3 cm

The corollary  below is a significant improvement from \cite{MS} of
our understanding of $T(\S )/Comm_{per}(\pi_1(T))$. Namely, we
showed in \cite{MS} that $T(\S )/Comm_{per}(\pi_1(T))$ is
non-Hausdorff by showing that orbits under $PSL_2(\mathbb{Z})$ of
marked hyperbolic metrics on $\S$ which are not lifts of hyperbolic
metrics on finite surfaces have accumulation points in $T(\S )$. In
this paper we start with the basepoint in $T(\S )$, i.e. a marked
hyperbolic metric from the Modular torus, and find an explicit
sequence of elements in $Comm_{per}(\pi_1(T))$ such that the image
of the basepoint under these elements accumulates onto itself.
Moreover, we establish that the orbit of the basepoint has closure
strictly larger that the orbit itself.

\vskip .3 cm

\paragraph{\bf Corollary 2} {\it The closure in the Teichm\"uller metric
of the orbit (under the base leaf preserving mapping class group
$Comm_{per}(\pi_1(T))$) of the basepoint in $T(\S )$ is strictly
larger than the orbit. Moreover, the closure of this orbit is
uncountable.}

\vskip .3 cm

The above Corollary is proved using the Baire category theorem and
Theorem 3.3 (see Section 3). However, we are also able to find an
explicit sequence in $Comm_{per}(\pi_1(T))$ whose limit point in
$T(\S )$ is not an element of $Comm_{per}(\pi_1(T))$ (see Corollary
4.2 in Section 4).

\section{The Farey tessellation}

We define the Farey tessellation $\mathcal{F}$ of the unit disk $\D$
as follows (see Figure 1). Let $\Delta_0$ be the ideal triangle in
$\D$ with vertices $-1$, $1$ and $i$. We invert $\Delta_0$ by
applying the three hyperbolic involutions, each of the three
preserves setwise one boundary side of $\Delta_0$ (but it changes
the orientation on the corresponding geodesic). By this, we obtain
three more ideal triangles each sharing one boundary side with
$\Delta_0$. We continue the inversions with respect to the new
triangles indefinitely. As a result, we obtain a locally finite
ideal triangulation of $\D$ called {\it the Farey tessellation}
$\mathcal{F}$. The set of the vertices in $S^1$ of the ideal
triangles from $\mathcal{F}$ is denoted by $\bar{\mathbb{Q}}$.  A
hyperbolic geodesic that is a side of a triangle from  $\mathcal{F}$
is also called an {\it edge} in $\mathcal{F}$. Denote by $l_0$ the
edge with the endpoints $-1$ and $1$, and fix an orientation on
$l_0$ such that $-1$ is the initial point and $1$ is the terminal
point. We call this edge the {\it distinguished oriented edge} of
$\mathcal{F}$. Also, denote by $l_1$ the oriented edge of
$\mathcal{F}$ with the endpoints $1$ and $i$ (and in that order).

\vskip .3 cm

Let $f:S^1 \to S^1$ be a homeomorphism. Then $f(\mathcal{F})$ is a
well defined ideal triangulation of  $\D$. We say that $\mathcal{F}$
is invariant under $f$ if $f(\mathcal{F})=\mathcal{F}$ as the ideal
triangulations. The Farey tessellation $\mathcal{F}$ is invariant
under the action of the group $PSL_2(\mathbb{Z})$. If a
homeomorphism of $S^1$ preserves $\mathcal{F}$, then it is
necessarily in $PSL_2(\mathbb{Z})$. This easy but important
observation  was proved in \cite{PS}.

\begin{figure}
\centering
\includegraphics[scale=0.6]{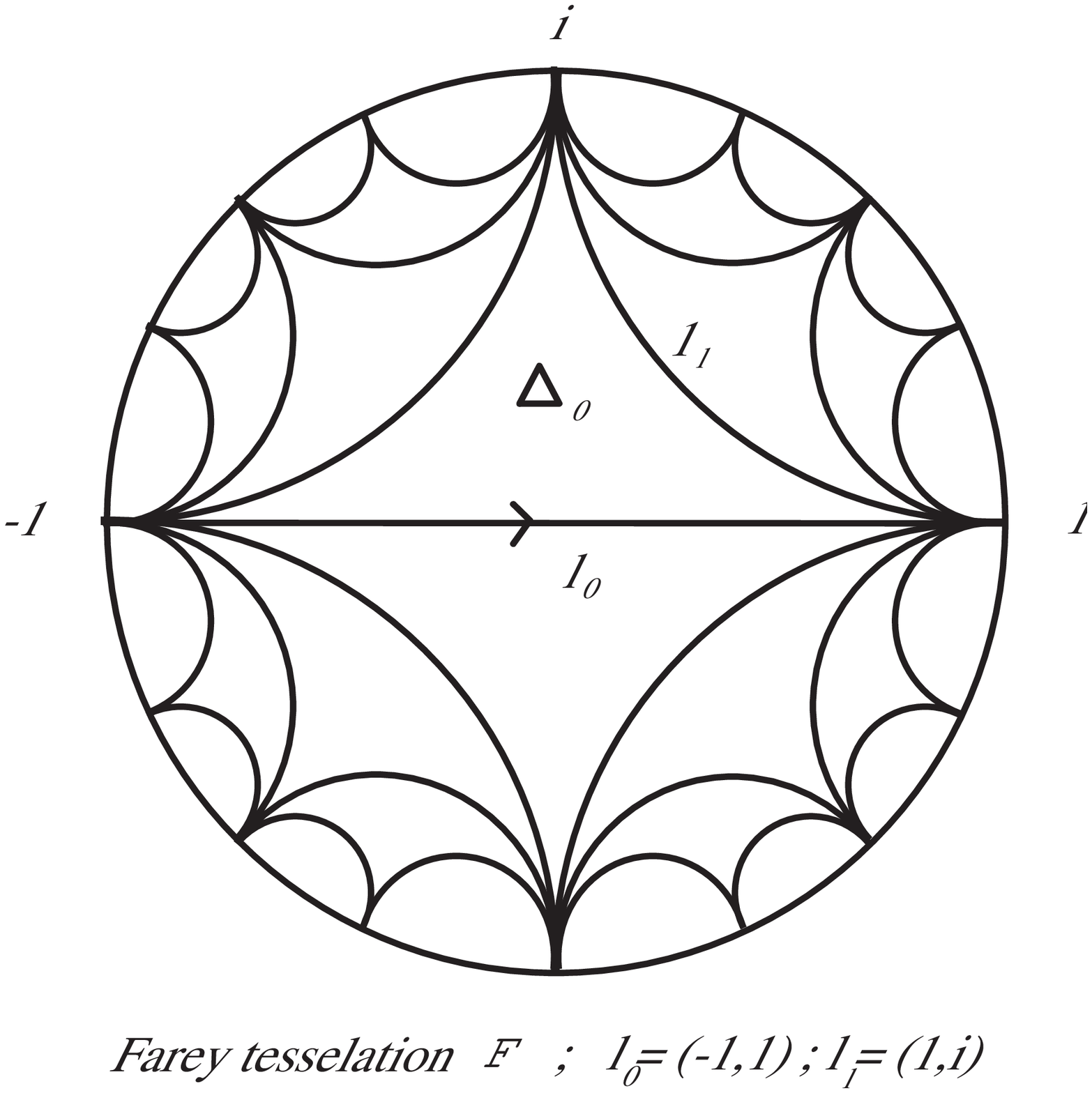}
\caption{} \label{}
\end{figure}

\vskip .3 cm

Consider two arbitrary locally finite ideal triangulations
$\mathcal{F}^1$ and $\mathcal{F}^2$ of $\D$. Fix two oriented edges
$e_1$ and $e_2$ from $\mathcal{F}^1$ and $\mathcal{F}^2$,
respectively. Then $e_1$ and $e_2$ are called the {\it distinguished
oriented edges} of tessellations $\mathcal{F}^1$ and
$\mathcal{F}^2$, respectively. There exists a unique homeomorphism
$h:S^1\to S^1$ which maps $\mathcal{F}^1$ onto $\mathcal{F}^2$ such
that $e_1$ is mapped onto $e_2$ in the orientation preserving manner
\cite{P}. We call such $f$ the {\it characteristic map} of
$\mathcal{F}^1$ and $\mathcal{F}^2$. (The characteristic maps
between the Farey tessellation and an arbitrary tessellation of $\D$
were used in \cite{P} to study the space of homeomorphisms of $S^1$.
In this paper, we use the notion of a characteristic map in a
slightly broader sense that its domain is not only the Farey
tessellation, but we allow an arbitrary tessellation.)

\vskip .3 cm

We recall the construction of  Whitehead homeomorphisms of $S^1$
(the construction below has been developed in \cite{PS}). Throughout
this paper $G_0<PSL_2(\mathbb{Z})$ denotes the finite index subgroup
such that $\D/G_0$ is the Modular torus $T_0$. An ideal
triangulation of $\D$ is said to be an {\it invariant tessellation}
if it is invariant under the action of a finite index subgroup
$K<G_0$. Equivalently, an invariant tessellation is an ideal
triangulation of $\D$ that is the lift of a finite, ideal
triangulation of some finite Riemann surface that covers  $T_0$. In
particular, the Farey tessellation is an invariant tessellation.

We use the following result:

\vskip .3 cm

\paragraph{\bf Theorem 2.1} (\cite{PS}) {\it
 Let $\mathcal{F}^1$ and $\mathcal{F}^2$ be two invariant tessellations
with the distinguished oriented edges $e_1$ and $e_2$, respectively.
The characteristic map of $\mathcal{F}^1$ onto $\mathcal{F}^2$ which
sends the distinguished oriented edge $e_1$ onto the distinguished
oriented edge $e_2$ conjugates a finite index subgroup of $G_0$ onto
another (possibly different) finite index subgroup of $G_0$.}

\vskip .3 cm

Let $\mathcal{F}^1$ be an invariant tessellation (with the
distinguished oriented edge $e_1$) that is invariant under the
action of a finite index subgroup $K<G_0$ and let $e$ be an edge of
$\mathcal{F}^1$. Let $K_1<K$ be a finite index subgroup. For the
simplicity of the exposition, we assume that the distinguished
oriented edge $e_0$ does not belong to the orbit $K_1\{ e\}$ of the
edge $e$ of $\mathcal{F}^1$.

\vskip .1 cm

\paragraph{\bf Definition 2.2} A {\it Whitehead move} on $\mathcal{F}_1$
along the orbit $K_1\{ e\}$ is the operation of replacing the orbit
of edges $K_1\{ e\}$ by the new orbit of edges $K_1\{ f\}$, where
$f$ is other diagonal of the ideal quadrilateral in $(\D \setminus
\mathcal{F}_1) \cup  \{ e\}$ (see Figure 2). As the result of this
operation we obtain a new ideal triangulation of the unit disc $\D$
that is in fact an invariant tessellation. This new tessellation is
denoted by $\mathcal{F}^1_{K_1, e}$, it is invariant under the
action of the group $K_1$ and its distinguished oriented edge is
$e_0$.

\vskip .3 cm

Consider the homeomorphisms $h$ of $S^1$ which fixes $e_0$ and which
maps $\mathcal{F}^1$ onto $\mathcal{F}^1_{K_1,e}$. By Theorem 2.1,
$h$ conjugates a finite index subgroup of $G_0$ onto another finite
index subgroup of $G_0$.

\vskip .1 cm

\paragraph{\bf Definition 2.3} (\cite{PS}) The {\it Whitehead
homeomorphism} corresponding to the Whitehead move along the orbit
$K_1\{ e\}$ is the characteristic map $h:S^1\to S^1$ of
$\mathcal{F}^1$ and $\mathcal{F}^1_{K_1, e}$ that fixes the common
distinguished oriented edge $e_0$ of $\mathcal{F}^1$ and
$\mathcal{F}^1_{K_1, e}$. It follows directly that the Whitehead
homeomorphism $h$ depends only on $\mathcal{F}^1$ and
$\mathcal{F}^1_{K_1, e}$, and we already noted that it conjugates a
finite index subgroup of $G_0$ onto another (possibly different)
finite index subgroup of $G_0$.

\begin{figure}
\centering
\includegraphics[scale=0.5]{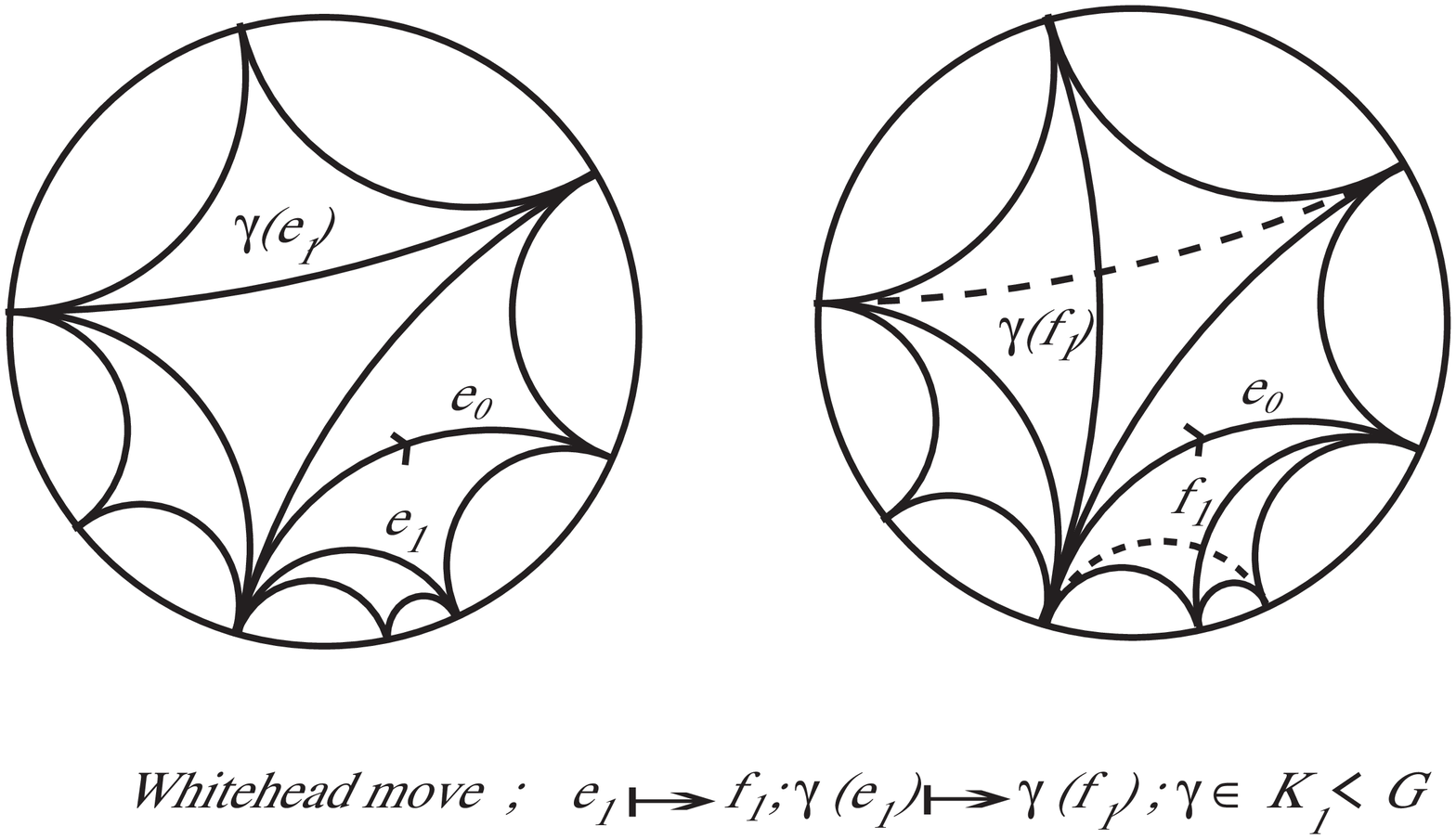}
\caption{} \label{}
\end{figure}

\vskip .3 cm

\section{Mapping Classes with Small Dilatations}

The following lemma is the first step in finding quasiconformal maps
of $\D$ which conjugate two finite index subgroups of
$PSL_2(\mathbb{Z})$ that are not conformally conjugated to each
other. We first show that the barycentric extensions \cite{DE} of
the Whitehead homeomorphisms (which we have defined in the previous
section) have dilatations essentially supported in a neighborhood of
the diagonal exchange for the corresponding Whitehead move.

\vskip .3 cm We will use the following notation below.

\paragraph{\bf Definition 3.1} Let $F:\D \to \D$ be a quasiconformal map
and let $N\in\mathbb{N}$. Then

$$
V(F,N):=\{ z\in\D:\ |\mu(F)(z)|\geq\frac{1}{N}\} .
$$

\vskip .3 cm

Let $f$ be a homeomorphism of the circle. By $E(f):\D \to \D$ we
always denote the barycentric extension of $f$ (see \cite{DE}).

\vskip .3 cm

Recall that $\mathcal{F}$ is the Farey tessellation with the
distinguished oriented edge $l_0$ (which is an oriented geodesic
with endpoints $-1$ and $1$). We keep the notation $l_1$ for the
edge of $\mathcal{F}$ whose endpoints are $1$ and $i$. Let $A\in
PSL_2(\mathbb{Q})$ be a hyperbolic translation with the oriented
axis $l_0$.  Let $\mathcal{F}^A$ denote the image $A(\mathcal{F})$
of the Farey tessellation $\mathcal{F}$ under $A$. Then
$\mathcal{F}^A$ is invariant under the group $AG_{0}A^{-1}$.  Define
$G_A:=G_0\cap AG_{0}A^{-1}$. Since $A$ is in the commensurator of
$PSL_2(\mathbb{Z})$ and since intersections of finitely many finite
index subgroups of $PSL_2(\mathbb{Z})$ is a finite index subgroup of
$PSL_2(\mathbb{Z})$, we conclude that the group $G_A$ is a subgroup
of finite index in $PSL_2(\mathbb{Z})$. It follows that
$\mathcal{F}^A$ is an invariant tessellation of $\D$ which is
invariant under the finite index subgroup $G_A<G_0$ (note that the
relation $G_A<G_0$ follows from the definition of $G_A$).

\vskip .3 cm

\paragraph{\bf Lemma 3.2} {\it Let $\mathcal{F}^A$ be an invariant
tessellation of $\D$ which is the image of the Farey tessellation
$\mathcal{F}$ under a hyperbolic translation $A\in
PLS_2(\mathbb{Q})$ with the oriented axis $l_0$. Let
$\mathcal{F}_{G,A(l_1)}^A$ be the image of $\mathcal{F}^A$ under the
Whitehead move along the orbit $G\{ A(l_1)\}$, where $G<G_A$ is any
subgroup of finite index. Let $f_A$ be the Whitehead homeomorphism
which maps $\mathcal{F}_{G,A(l_1)}^A$ onto $\mathcal{F}^A$ fixing
the common distinguished oriented edge $l_0$ and let $E(f_A)$ be its
barycentric extension. Then, for each $N\in\mathbb{N}$ there exists
$K_N=K_N(z_0,\mathcal{F}^A)>0$ such that $V(E(f_A),N)$ is a subset
of the $K_N$-neighborhood of the orbit $G\{ z_0\}$, where $z_0\in
l_0$ is an arbitrary point. The constant $K_N$ is independent of
$G$.}

\vskip .3 cm

\paragraph{\bf Remark } According to the definition of the
Whitehead homeomorphisms, the characteristic map between
$\mathcal{F}^A$ and $\mathcal{F}^A_{G,A(l_1)}$ is the Whitehead
homeomorphism, and $f_A$ is its inverse. However, the Whitehead move
on $\mathcal{F}^A_{G,A(l_1)}$ along the orbit $G\{ A(l_1')\}$ gives
$\mathcal{F}^A$, where $l_1'$ is the other diagonal of the ideal
quadrilateral in $(\D \setminus \mathcal{F})\cup \{ l_1\}$.
Therefore, $f_A$ is also a Whitehead homeomorphism corresponding to
this ``inverse'' Whitehead move. Although the notation $f_A$ does
not suggest that the map $f_A$ depends on the group $G<G_A$, it is
important to remember that it does. It will always be clear from the
context what is the corresponding group $G$.

\vskip .3 cm

\paragraph{\bf Remark }
The above lemma includes the possibility that $A=id$. In this case
the barycentric extension of the Whitehead homeomorphism $f_{id}$
between $\mathcal{F}_{G,l_1}$ and the Farey tessellation
$\mathcal{F}$ is supported on the $G$ orbit of a neighborhood of
$z_0\in l_0$.

\vskip .3 cm

\paragraph{\bf Remark} The inverse $E(f_A)^{-1}$ of the barycentric
extension of the Whitehead homeomorphism $f_A$ which maps
$\mathcal{F}^A$ onto $\mathcal{F}^A_{G,A(l_1)}$ has Beltrami
dilatation ``essentially'' supported in
$K_N'=K_N'(K_N)$-neighborhood of $E(f_A)(G\{ z_0\} )$. If $A=id$
then $E(f_{id})(G\{ z_0\} )=H\{ E(f_{id})(z_0)\}$, where
$H<PSL_2(\mathbb{Z})$ is conjugated to $G$ by $f_{id}$ (see Theorem
2.1 and its proof in \cite{PS}).

\vskip .3 cm

\paragraph{\bf Proof} Fix $N \in \mathbb{N}$. The proof is by
contradiction. That is we assume that there exists a sequence of
subgroups $G_n<G_A$ of finite index (every element in every group
$G_n$ necessarily preserves $\mathcal{F}^A$) and a sequence of
points $w_n\in\D$ such that
$$
dist(w_n,G_n\{ z_0\} )\to\infty
$$
and
$$
|\mu(E(f_n))(w_n)|\geq \frac{1}{N},
$$
where $f_n$ is the Whitehead homeomorphism which maps
$\mathcal{F}^A_{G_n,A(l_1)}$ onto $\mathcal{F}^A$ fixing the common
distinguished oriented edge $l_0$.

\vskip .3 cm

Let $\mathcal{Z}_0:=G_A\{ z_0\}$ be the full orbit of $z_0$ under
$G_A$. After passing onto a subsequence if necessary, there are two
cases that we have to consider:

\begin{enumerate}
\item  There exists $C>0$ so that $dist(w_n,\mathcal{Z}_0)\leq C$ for all
$n \in \mathbb{N}$.

\item We have that $dist(w_n,\mathcal{Z}_0)\to\infty$ as $n\to\infty$.

\end{enumerate}

We denote by $\mathcal{F}_n^A:=\mathcal{F}_{G_n,A(l_1)}^A$ the image
of the invariant tessellation $\mathcal{F}^A$ under the Whitehead
move along the orbit $G_n\{ A(l_1)\}$ of $A(l_1)$.

\vskip .3 cm

We first settle the first case, that is we assume that
$dist(w_n,\mathcal{Z}_0)\leq C$ for all $n$. Since $l_1$ is within
the bounded distance from $z_0$ it follows that $w_n$ is within the
bounded distance from $G_A\{ l_1\}$. From the assumptions that
$dist(w_n,G_n\{ z_0\} )\to\infty$ as $n\to\infty$, and that
$dist(w_n,\mathcal{Z}_0)\leq C$ for all $n$, we get that
$dist(w_n,G_n\{ l_1\} )\to\infty$ as $n\to\infty$.

\vskip .3cm

Recall that $f_n$ is the Whitehead homeomorphism which maps
$\mathcal{F}_n^A$ onto $\mathcal{F}^A$ and which fixes the common
distinguished oriented edge $l_0$. The barycentric extension of
$f_n$ is denoted by $E(f_n)$.

\vskip .3 cm

Choose $\gamma_n \in G_A$ such that $dist(w_n,\gamma_n(z_0))\leq C$
for all $n\in \mathbb{N}$. Since $\mathcal{F}^A$ is invariant under
$G_A$, there exists a fundamental polygon for $G_A$ which is a union
of finitely many adjacent triangles from $\mathcal{F}^A$. Moreover,
we can choose such a fundamental polygon $\omega$ with the following
properties

\begin{enumerate}
\item The boundary of $\omega$ contains the distinguished oriented
edge $l_0$.

\item The polygon $\omega$ is to the left of $l_0$.

\item We have $l_1\subset\omega^{\circ}$, where $\omega^{\circ}$ is the
interior of $\omega$.
\end{enumerate}

The union of translates of $\omega$ under the group $G_A$ tiles the
unit disk $\D$.

\vskip .3 cm It is important to note that the tessellations
$\mathcal{F}^A$ and $\mathcal{F}_n^A$ agree on the orbit
$(G_A\setminus G_n)\{\omega \}$ of the fundamental polygon $\omega$
(they differ inside the orbit $G_n\{\omega \}$). Let $\alpha_n\in
G_A$ be such so that $w_n\in\alpha_n(\omega )$. Also, let
$T_n\subset\alpha_n(\omega )$ be a triangle in $\mathcal{F}_n^A$
which contains $w_n$. Then the triangles $\alpha_n^{-1}(T_n)$ are
contained in $\omega$ for each $n$. After passing onto  a
subsequence if necessary we may assume that $\alpha_n^{-1}(T_n)$ is
the same triangle $T$ in $\omega$ for each $n$.

\vskip .3 cm

Since $dist(w_n,G_n\{ l_1\} )\to\infty$ as $n\to\infty$ and
$dist(w_n,G_A\{ z_0\})\leq C$ for all $n$, we conclude  that the
tessellations $\alpha_n^{-1}(\mathcal{F}_n^A)$ and $\mathcal{F}^A$
agree on the edges intersecting a hyperbolic disk with the center
$z_0$ and the hyperbolic radius $r_n$, where $r_n\to\infty$ as
$n\to\infty$. In particular, the triangle $T$ is in $\mathcal{F}^A$.
We already know that the triangle $T'_n=f_n(T_n)$ is in
$\mathcal{F}^A$ because $f_n$ maps $\mathcal{F}_n^A$ onto
$\mathcal{F}^A$. Therefore, there exists a unique $\beta_n \in
A\circ PSL_2(\mathbb{Z})\circ A^{-1}$ such that $\beta_n(T_n')=T$
and such that $\beta_n \circ f_n \circ \alpha_n$ fixes each vertex
of $T$ (the fact that $PSL_2(\mathbb{Z})$ is transitive on the
oriented edges of the Farey tessellation $\mathcal{F}$ implies that
$A\circ PSL_2(\mathbb{Z})\circ A^{-1}$ is transitive on the oriented
edges of the invariant tessellation $\mathcal{F}^A$ which implies
the existence of such $\beta_n$).

\vskip .3 cm

The circle homeomorphism  $\beta_n\circ f_n\circ \alpha_n$ maps
$\alpha_n^{-1}(\mathcal{F}_n^A)$ onto
$\beta_n(\mathcal{F}^A)=\mathcal{F}^A$. Its barycentric extension is
$\beta_n\circ E(f_n)\circ \alpha_n$ (the barycentric extension is
conformally natural, see \cite{DE}). As we have already shown, given
any neighborhood of the origin in the unit disc, we can find $n \in
\mathbb{N}$, so that the tessellations
$\alpha_n^{-1}(\mathcal{F}_n^A)$ and $\mathcal{F}^A$ agree on that
neighborhood. Since $\beta_n\circ f_n\circ \alpha_n$ fixes every
vertex of the triangle $T$ it follows that $\beta_n\circ f_n\circ
\alpha_n \to id$ on the circle as $n \to \infty$. This implies that
the Beltrami dilatation $\mu (\beta_n\circ E(f_{n})\circ \alpha_n)$
converges to zero uniformly on compact subsets of $\D$. Since $w_n$
is on the bounded distance from $\mathcal{Z}_0$ it follows that
$\alpha_n^{-1}(w_n)$ is in a compact subset of $\D$. This implies
that
\begin{equation}
\label{one} |\mu (\beta_n\circ E(f_n)\circ
\alpha_n)(\alpha_n^{-1}(w_n))|\to 0
\end{equation}
as $n\to\infty$. Since $|\mu (\beta_n\circ E(f_{n})\circ
\alpha_n)|=|\mu (E(f_n))\circ\alpha_n|$, we derive a contradiction
to the assumption that $|\mu(E(f_n))(w_n)|\geq\frac{1}{N}$ for all
$n\in\mathbb{N}$.

\vskip .3 cm

It remains to consider the case when
$dist(w_n,\mathcal{Z}_0)\to\infty$ as $n\to\infty$. We keep the
notation $f_n$ for the Whitehead homeomorphism which maps
$\mathcal{F}_n^A:=\mathcal{F}_{G_n,A(l_1)}^A$ onto $\mathcal{F}^A$
and which fixes $l_0$.

\vskip 1cm

\paragraph{\bf Remark } Note that the condition
$dist(w_n,\mathcal{Z}_0)\to\infty$, as $n\to\infty$, means that the
projection of the sequence $w_n$ onto the surfaces obtained as the
quotient of the unit disc by finite index subgroups $G_n$ of $G_0$
converges to the punctures in the boundary of that surfaces. The
fact (that we prove in detail below) that the Beltrami dilatation of
$E(f_n)$ tends to zero along the sequence $w_n$ is actually a
corollary of the fact that the circle homeomorphism $f_n$ is
differentiable at every ``rational'' point on the circle (these are
the fixed points of parabolic transformations from $G_0$).

\vskip .3 cm

We fix a fundamental polygon $\omega$ for $G_A$ as above. That is,
$\omega$ is the union of adjacent triangles in $\mathcal{F}^A$ such
that
\begin{enumerate}
\item $l_0$ is on the boundary of the fundamental polygon $\omega$
\item $\omega$ is to the left of $l_0$
\item $l_1\subset\omega^{\circ}$, where $\omega^{\circ}$ is the interior
of $\omega$
\end{enumerate}

Let $T_n$ be a triangle from $\mathcal{F}_n^A$ which contains $w_n$
and let $T_n'=f_n(T_n)$ be the image triangle in $\mathcal{F}^A$ as
before. Let $\alpha_n\in G_A$ be such that
$T_n\subset\alpha_n(\omega )$ and let $\beta_n\in G_A$ be such that
$\beta_n(T_n')\subset\omega$. After passing onto a subsequences if
necessary, we can assume that $\alpha_n^{-1}(T_n)$ and
$\beta_n(T_n')$ are fixed triangles $T$ and $T'$ in $\omega$,
respectively. We note that $T'$ is a triangle in $\mathcal{F}^A$. On
the other hand, $T$ is in $\alpha_n^{-1}(\mathcal{F}^A_n)$ for each
$n$, which implies that that $T$ is either in $\mathcal{F}^A$ (if
$\alpha_n\in G_A\setminus G_n$ or if $\alpha_n\in G_n$ and
$T_n\subset\mathcal{F}^A\setminus \mathcal{F}^A_n$) or in
$\mathcal{F}^A_n$ (if $\alpha_n\in G_n$ and
$T_n\subset\mathcal{F}^A_n\setminus\mathcal{F}^A$).

\vskip .3 cm

The map $\widetilde{f}_n:=\beta_n\circ f_n\circ \alpha_n^{-1}$ maps
$T$ onto $T'$. After passing onto a subsequence if necessary, the
points $\alpha_n^{-1}(w_n)\in T$ converge to a single ideal vertex
$y$ of $T$. Let $y':=\widetilde{f}_n(y)\in T'$ which is a fixed
point for all $n$ after possibly passing onto  a subsequence if
necessary. Let $l$ be a boundary side of $T$ with $y$ its ideal
endpoint such that $l\in\mathcal{F}^A$ (at least one of the two
boundary sides of $T$ with $y$ their ideal endpoint is in
$\mathcal{F}^A$ because the tessellation
$\alpha_n^{-1}(\mathcal{F}^A_n)$ is obtained by a Whitehead move on
$\mathcal{F}^A$ along the orbit $\alpha_n^{-1}G_n\{ A(l_1)\}$ which
implies that no two adjacent geodesics can be changed by the
definition of a Whitehead move). Let
$l'=\widetilde{f}_n(l)\in\mathcal{F}^A$ be a boundary side of $T'$
with an ideal endpoint $y'$ (since $\widetilde{f}_n(T)$ is a fixed
triangle $T'$ for each $n$, then after passing onto a subsequence if
necessary the side $\widetilde{f}_n(l)$ is the same boundary side
$l'$ of $T'$).

\vskip .3 cm

Let $\gamma\in G_A$ be a primitive parabolic element which fixes $y$
and let $\gamma'\in G_A$ be a primitive parabolic element in $G_A$
which fixes $y'$. Then the set of edges in $\mathcal{F}^A$ with one
ideal endpoint $y$ is invariant under the action of $\gamma$ and a
fundamental set for the action of a cyclic group $<\gamma >$
generated by $\gamma$ consists of finitely many adjacent geodesics
of $\mathcal{F}^A$ with one endpoint $y$. Similarly, the set of
edges in $\mathcal{F}^A$ with one ideal endpoint $y'$ is invariant
under the action of $\gamma'$ and a fundamental set for the action
of $<\gamma'>$ consists of finitely many adjacent geodesics of
$\mathcal{F}^A$ with one endpoint $y'$.

\vskip .3 cm

The group $G_n$ is a finite index subgroup of $G_A$ and it follows
that $\alpha_n^{-1}G_n\alpha_n$ is also a finite index subgroup of
$G_A$. Therefore, the isotropy subgroup of $y$ in
$\alpha_n^{-1}G_n\alpha_n$ is of finite index in the isotropy group
$<\gamma >$ of $y$ in $G_A$. Thus a generator $\gamma_n\in
\alpha_n^{-1}G_n\alpha_n$ of the isotropy group of $y$ is equal to a
finite, non-zero, integer power of $\gamma$. After possibly
replacing $\gamma_n$ by its inverse if necessary, we have
$\gamma_n=\gamma^{p_n}$, for some $p_n\in\mathbb{N}$. A fundamental
set for the action of $<\gamma_n>$ on the geodesics of
$\mathcal{F}^A$ with one endpoint $y$ is obtained by taking $p_n$
translates by $\gamma$ of a fixed fundamental set for $<\gamma>$.

\vskip .3 cm

Note that the tessellation $\alpha_n^{-1}(\mathcal{F}^A_n)$ is
obtained by a Whitehead move on $\mathcal{F}^A$ along the orbit
$\alpha_n^{-1}G_n\alpha_n\{ A(l_1)\}$. Let $k_y$ be the number of
geodesics in a fundamental set for the action of $<\gamma >$ on the
set of edges of $\mathcal{F}^A$ which have one endpoint $y$. If $l$
is a fixed edge of $\mathcal{F}^A$ with one endpoint $y$ then $k_y$
is the number of edges in $\mathcal{F}^A$ with one endpoint $y$
which lie in between $l$ and $\gamma (l)$, where we count $l$ but do
not count $\gamma (l)$. Let $k_{y'}$ be the number of geodesics in a
fundamental set for the action of $<\gamma'>$ on the edges of
$\mathcal{F}^A$ with one endpoint $y'$. Equivalently, $k_{y'}$ is
the number of geodesics in $\mathcal{F}^A$ with one ideal endpoint
$y'$ in between $l'$ and $\gamma' (l')$ counting $l'$ but not
counting $\gamma'(l')$, where $l'$ is a fixed edge of
$\mathcal{F}^A$ with one endpoint $y'$. The number of geodesics in a
fundamental set for the action of $<\gamma_n'>$ on the set of edges
of $\mathcal{F}^A$ with one endpoint $y$ is $k_yp_n$.

\vskip .3 cm

Recall that $\alpha_n^{-1}(\mathcal{F}^A_n)$ is obtained by a
Whitehead move on $\mathcal{F}^A$ along the orbit
$\alpha_n^{-1}G_n\alpha_n\{ A(l_1)\}$. The Whitehead move can either
add edges at $y$, erase edges at $y$, or do not change edges at $y$.
We further assume that the choice of the edge $l$ in $\mathcal{F}^A$
with one endpoint $y$ is such that the Whitehead move does not erase
$l$. The number of geodesics in $\alpha_n^{-1}(\mathcal{F}^A_n)$
with an ideal endpoint at $y$ in between $l$ and $\gamma_n(l)$
(including $l$ but not including $\gamma_n(l)$) is $k_yp_n+a$, where
$a=0$ if the Whitehead move does not change any edge at $y$, $a=1$
if the Whitehead move adds edges at $y$ or $a=-1$ if the Whitehead
move erases edges at $y$. (In the top part of Figure 3, we
illustrate the case when the Whitehead move adds geodesics at $y$;
$k_y=3$; $k_{y'}=2$.)

\begin{figure}
\centering
\includegraphics[scale=0.4]{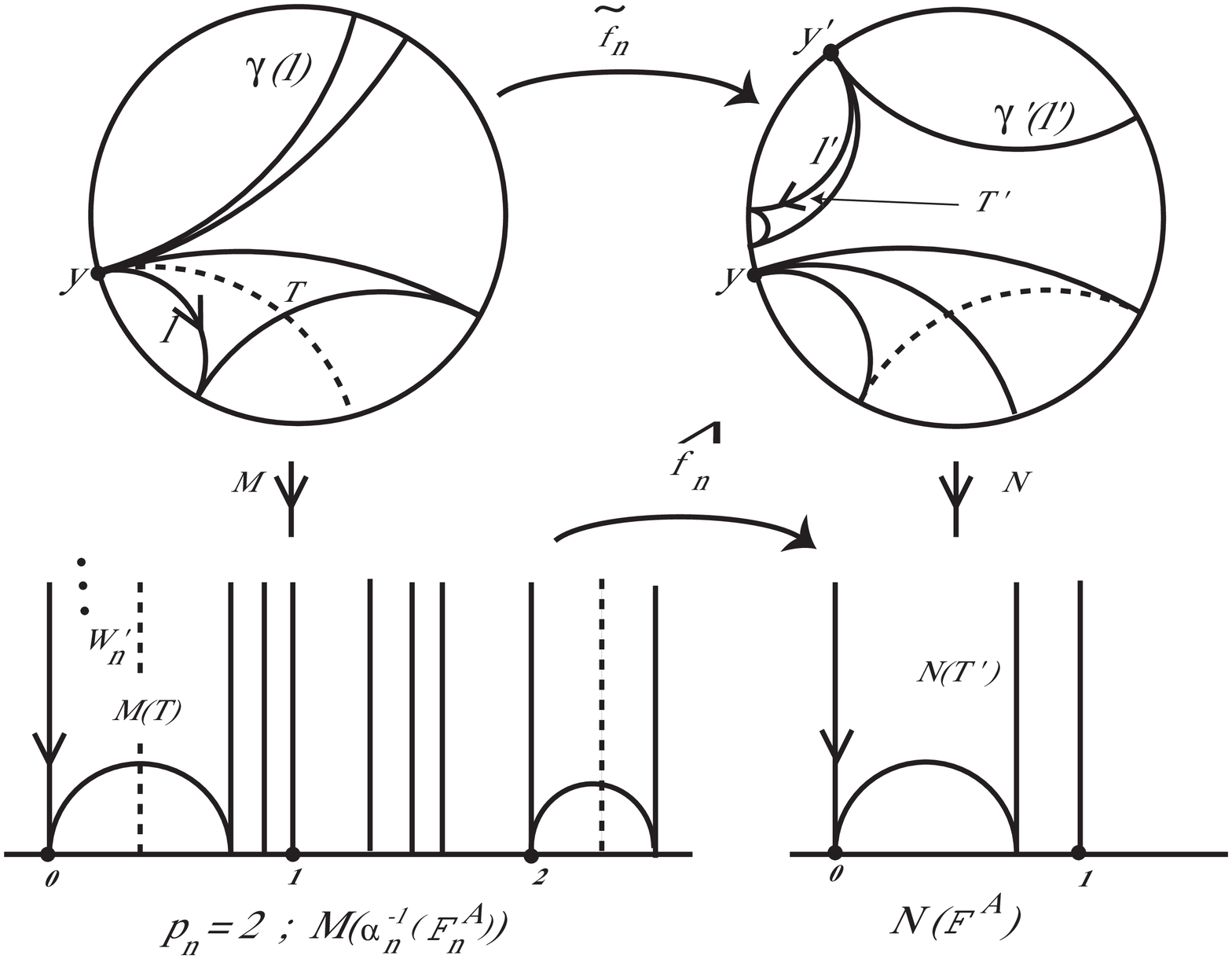}
\caption{} \label{}
\end{figure}

\vskip .3 cm

Let $M:\D\to\mathbb{H}$ be a M\"obius map which sends $y$ to
$\infty$, $l$ to a geodesic with endpoints $0$ and $\infty$, and
$\gamma (l)$ to a geodesic with endpoints $1$ and $\infty$. Let
$N:\D\to\mathbb{H}$ be a M\"obius map which sends $y'$ to $\infty$,
$l'$ to a geodesic with endpoints $0$ and $\infty$, and $\gamma'
(l')$ to a geodesic with endpoints $1$ and $\infty$. Define
$\widehat{f}_n:=N\circ \widetilde{f}_n\circ M^{-1}$. Let
$w_n':=M(\alpha_n^{-1}(w_n))\in\mathbb{H}$ (see Figure 3). Then
$w_n\to\infty$ as $n\to\infty$ inside the triangle $M(T)$. Namely,
$b_n:=Im(w_n')\to\infty$ as $n\to\infty$ and $0\leq Re(w_n')<1$ for
all $n\in\mathbb{N}$. This implies that $\frac{1}{b_n}w_n'$ stays in
a compact subset of $\mathbb{H}$.

\vskip .3 cm

We consider the pointwise limit of
$\frac{1}{b_n}\widehat{f}_n(b_nx)$ as $n\to\infty$ for all
$x\in\mathbb{R}$. Our goal is to show that it is a linear map. There
are two possibilities (after passing onto a subsequence if
necessary)
\begin{enumerate}
\item $p_n\to \infty$ as $n\to\infty$
\item $p_n=p$ is fixed, for all $n\in\mathbb{N}$
\end{enumerate}

\vskip .3 cm

We assume that $p_n\to\infty$ as $n\to\infty$. We obtain an upper
bound
\begin{equation}
\label{upper} \widehat{f}_n(b_nx)< \widehat{f}_n([b_nx]+1)\leq
([b_nx]+1)k_y\frac{1}{k_{y'}}+
\Big{(}\frac{[b_nx]+1}{p_n}+1\Big{)}\frac{1}{k_{y'}}+2,
\end{equation}
where $[b_nx]$ is the integer part of $b_nx$. The first inequality
in (\ref{upper}) follows because $b_nx<[b_nx]+1$ and $\widehat{f}_n$
is an increasing function. The second inequality in (\ref{upper}) is
obtained as follows. By the choice of $M$ and $N$ above, we have
$(M\circ \gamma\circ M^{-1})(z)=z+1$, $(M\circ \gamma_n\circ
M^{-1})(z)=z+p_n$ and $(N\circ \gamma'\circ N^{-1})(z)=z+1$ (see
Figure 3). In between $0$ and $[b_nx]+1$ there is
$[\frac{[b_nx]+1}{p_n}]$ of adjacent intervals of length $p_n$. For
each interval of length $p_n$, the number of geodesics in
$M(\alpha_n^{-1}(\mathcal{F}^A_n))$ with one endpoint $\infty$ and
the other endpoint in the interval is at most the number of
geodesics in $M(\mathcal{F}^A)$ with one endpoint $\infty$ and the
other point in the interval plus one extra geodesic (because the
Whitehead move adds at most one geodesic in such an interval). The
map $\widehat{f}_n$ fixes $0$ and $\infty$, and it maps the
geodesics of $M(\alpha_n^{-1}(\mathcal{F}^A_n))$ onto the geodesics
of $N(\mathcal{F}^A)$. Therefore, we need to estimate the number of
geodesics in $M(\alpha_n^{-1}(\mathcal{F}^A_n))$ with one endpoint
$\infty$ and the other endpoint in the interval $[0,[b_nx]+1]$. The
second inequality in (\ref{upper}) is obtained by noting that
$([b_nx]+1)k_y$ is the number of geodesics in $M(\mathcal{F}^A)$
with one endpoint $\infty$ and the other point in the interval
$[0,[b_nx]+1]$ and that we add at most $\frac{[b_nx]+1}{p_n}+1$
geodesics to get the corresponding geodesics of
$M(\alpha_n^{-1}(\mathcal{F}^A_n))$. We need to divide the number of
geodesics by $k_{y'}$ because $N(\mathcal{F}^A)$ has $k_{y'}$
geodesics with one endpoint $\infty$ and the other endpoint in a
fixed interval of length $1$. Since the quantities $([b_nx]+1)k_y$
and $\frac{[b_nx]+1}{p_n}+1$ are not necessarily divisible with
$k_{y'}$, we add $2$ to ensure that we have an upper bound in
(\ref{upper}).

\vskip .3 cm

In a similar fashion, we obtain a lower bound
\begin{equation}
\label{lower} \widehat{f}_n(b_nx)\geq \widehat{f}_n([b_nx])\geq
[b_nx]k_y\frac{1}{k_{y'}}-\Big{(}\frac{[b_nx]}{p_n}+1\Big{)}\frac{1}{k_{y'}}-2.
\end{equation}

\vskip .3 cm

Since $p_n\to\infty$ and $\frac{[b_nx]}{b_n}\to x$ as $n\to\infty$,
the inequalities (\ref{upper}) and (\ref{lower}) imply that
$$
\frac{1}{b_n}\widehat{f}_n(b_nx)\to \frac{k_y}{k_{y'}}x,
$$
as $n\to\infty$. Thus $\frac{1}{b_n}\widehat{f}_n(b_nx)$ converges
pointwise to a linear map in the case when $p_n\to\infty$ as
$n\to\infty$.

\vskip .3 cm

We assume that $p_n=p$ for all $n\in\mathbb{N}$. Then we obtain the
following upper bound
\begin{equation}
\label{upper1} \widehat{f}_n(b_nx)<\widehat{f}_n([b_nx]+1)\leq
\Big{[}\frac{[b_nx]+1}{p}\Big{]}(k_yp+a)\frac{1}{k_{y'}}+k_yp\frac{1}{k_{y'}}.
\end{equation}
The second inequality in (\ref{upper1}) is obtained by noting that
there is $\Big{[}\frac{[b_nx]+1}{p}\Big{]}$ adjacent disjoint
intervals of length $p$ from $0$ to $[b_nx]+1$ each of which
contains endpoints of $k_yp+a$ geodesics of
$M(\alpha_n^{-1}(\mathcal{F}^A_n))$ with the other endpoint at
$\infty$. Since each interval of length $1$ contains $k_{y'}$
endpoints of geodesics of $N(\mathcal{F}^A)$ with the other endpoint
$\infty$, we obtain the first summand on the right of
(\ref{upper1}). We add $k_yp\frac{1}{k_{y'}}$ to the right of
(\ref{upper1}) because $\frac{[b_nx]+1}{p}$ might not be an integer
and, in this case, the interval $[[\frac{[b_nx]+1}{p}]p,[b_nx]+1]$
is of the length strictly smaller than $p$. Our upper estimate of
this part uses an interval of length $p$.

\vskip .3 cm

The following lower bound
\begin{equation}
\label{lower1} \widehat{f}_n(b_nx)\geq \widehat{f}_n([b_nx])\geq
\Big{[} \frac{[b_nx]}{p}\Big{]}(k_yp+a)\frac{1}{k_{y'}}
\end{equation}
is obtained similarly to the above upper bound.

\vskip .3 cm

The inequalities (\ref{upper1}) and (\ref{lower1}) together with the
facts that $[\frac{[b_nx]+1}{p}]/b_n\to \frac{1}{p}x$ and
$[\frac{[b_nx]}{p}]/b_n\to \frac{1}{p}x$ as $n\to\infty$ imply that
$$
\frac{1}{b_n} \widehat{f}_n(b_nx)\to
(\frac{k_y}{k_{y'}}+\frac{a}{pk_{y'}})x
$$
as $n\to\infty$. Thus $\frac{1}{b_n} \widehat{f}_n(b_nx)$ converges
to a linear map in the case $p_n=p$ as well.

\vskip .3 cm

We showed above that $\frac{1}{b_n}\widehat{f}_n( b_nx)$ converges
pointwise to a linear map in both cases which implies that $|\mu
(E(B_n^{-1}\circ \widehat{f}_n\circ B_n))|\to 0$ uniformly on
compact subsets, where $B_n(z):=b_nz$. Since $\frac{1}{b_n}w_n'$
stays in a compact subset of $\mathbb{H}$, we get (by the conformal
naturality of the barycentric extension \cite{DE}) that $|\mu
(E(B_n^{-1}\circ \widehat{f}_n\circ B_n))(\frac{1}{b_n}w_n')|=|\mu
(E(\widehat{f}_n))(w_n')|\to 0$ as $n\to\infty$. But this is in
contradiction with the starting assumption that $|\mu
(E(f_n))(w_n)|\geq\frac{1}{N}$ which proves the lemma. $\Box$

\vskip .3 cm

Recall that $A\in PSL_2(\mathbb{Q})$ is a hyperbolic translation
with the oriented axis $l_0$ and that $G_A=AG_0A^{-1}\cap G_0$. Then
$G_A$ is a finite index subgroup of $G_0$. We keep the notation
$\mathcal{F}$ for the Farey tessellation and the notation
$\mathcal{F}^A$ for the image of $\mathcal{F}$ under $A$. Then
$\mathcal{F}^A$ is a tessellation of $\D$ invariant under $G_A$.

\vskip .3 cm

If $G$ is a finite index subgroup of $G_0$, recall that
$\mathcal{F}_{G,l_1}$ is the image of $\mathcal{F}$ under the
Whitehead move along the orbit $G\{ l_1\}$. If $G$ is a finite index
subgroup of $G_A$, recall that $\mathcal{F}^A_{G,A(l_1)}$ is the
image of $\mathcal{F}^A$ under the Whitehead move along the orbit
$G\{ A(l_1)\}$.

\vskip .3 cm

We say that a sequence $f_n$ of quasisymmetric maps of $S^1$ {\it
converges in the Teichm\"uller metric} to a quasisymmmetric map $f$
of $S^1$ if there exists a sequence of quasiconformal extensions
$F_n:\D\to\D$ of $f_n$ and a quasiconformal extension $F:\D\to\D$ of
$f$ such that $\|\mu (F_n)-\mu (F)\|_{\infty}\to 0$ as $n\to\infty$.
Note that the Teichm\"uller metric on the space of quasisymmetric
maps of $S^1$ is a pseudometric. The Teichm\"uller metric projects
to a proper metric on the quotient of the space of quasisymmetric
maps of $S^1$ by the action of $PSL_2(\mathbb{R})$ (where the action
is given by the post-composition of quasisymmetric maps  with
$PSL_2(\mathbb{R})$).

\vskip .3 cm

We show below that the Whitehead homeomorphism from
$\mathcal{F}^A_{G,A(l_1)}$ to $\mathcal{F}^A$ followed by the
Whitehead homeomorphism from $\mathcal{F}$ to $\mathcal{F}_{G,l_1}$
converges to the identity in the Teichm\"uller metric as $A$
converges to the identity.

\vskip .3 cm

\paragraph{\bf Theorem 3.3} {\it Let $A\in PSL_2(\mathbb{Q})$ be a
hyperbolic translation with the oriented axis $l_0$. Let $G_A$ and
$\mathcal{F}^A$ be as above, and let $G$ be a finite index subgroup
of $G_A$. Let $f_{id}$ be the Whitehead homeomorphism fixing $l_0$
which maps $\mathcal{F}_{G,l_1}$ onto $\mathcal{F}$, and let $g_A$
be the Whitehead homeomorphism fixing $l_0$ which maps
$\mathcal{F}^A_{G,A(l_1)}$ onto $\mathcal{F}^A$. Then $$ g_A\circ
f_{id}^{-1}\to id$$ in the Teichm\"uller metric as $A\to id$. }

\vskip .3 cm

\paragraph{\bf Proof} We denote by $E(f_{id})$ and $E(g_A)$ the
barycentric extensions of $f_{id}$ and $g_A$, respectively. It is
enough to show that
$$
\|\mu (E(f_{id}))-\mu (E(g_A))\|_{\infty}\to 0
$$
as $A\to id$. (It is important to note that $\|\mu (E(f_{id}))-\mu
(E(g_A))\|_{\infty}$ is equal to $\sup_{z\in\D}|\mu
(E(f_{id}))(z)-\mu (E(g_A))(z)|$ because the barycentric extensions
of quasisymmetric maps are analytic diffeomorphisms which implies
that their Beltrami dilatations are continuous maps.)

\vskip .3 cm

Assume on the contrary that there exist sequences $w_n\in \D$,
$A_n\in PSL_2(\mathbb{Q})$ and $G_n<G_{A_n}$ such that $A_n$ is a
hyperbolic translation with the oriented axis $l_0$, $A_n\to id$ as
$n\to\infty$, $[G_{A_n}:G_n]<\infty$ and
\begin{equation}
\label{difference} |\mu (E(f_{id}))(w_n)-\mu (E(g_{A_n}))(w_n)|\geq
\frac{1}{N}
\end{equation}
for all $n\in \mathbb{N}$ and for a fixed $N\in \mathbb{N}$. This
implies that either $|\mu (E(f_{id}))(w_n)|$ or $|\mu
(E(g_{A_n}))(w_n)|$ is at least $1/N$. By Lemma 3.2, there exists
$K_N(z_0,\mathcal{F})>0$ such that $V(E(f_{id}))$ is a subset of the
$K_N(z_0,\mathcal{F})$-neighborhood of the orbit $G_n\{ z_0\}$.
Again by Lemma 3.2, there exists $K_N(z_0,\mathcal{F}^{A_n})>0$ such
that $V(E(g_{A_n}))$ is a subset of the
$K_N(z_0,\mathcal{F}^{A_n})$-neighborhood of the orbit $G_n\{
z_0\}$.

\vskip .3 cm

\paragraph{\bf Remark} Let $l_1'$ be the diagonal of the ideal
quadrilateral in $(\D\setminus\mathcal{F})\cup\{ l_1\}$ different
from $l_1$. We note that the Whitehead homeomorphism $f_{id}^{-1}$
from $\mathcal{F}$ to $\mathcal{F}_{G_n,l_1}$ does not necessarily
map the orbit $G_n\{ l_1\}$ in $\mathcal{F}$ onto the orbit $G_n\{
l_1'\}$; the Whitehead homeomorphism $g_{A_n}^{-1}$ from
$\mathcal{F}^{A_n}$ to $\mathcal{F}^{A_n}_{G_n, A(l_1)}$ does not
necessarily map the orbit $G_n\{ A_n(l_1)\}$ in $\mathcal{F}^{A_n}$
onto the orbit $G_n\{ A_n(l_1')\}$ in $\mathcal{F}^A_{G_n,A(l_1)}$.
On the other hand, $\mathcal{F}$ and $\mathcal{F}^{A_n}$ are both
obtained by infinite number of inversions in any of their triangles.
Thus, it is better to consider inverse Whitehead homeomorphisms
$f_{id}$ and $g_{A_n}$ because the image tessellations of
$\mathcal{F}_{G_n,l_1}$ and $\mathcal{F}^{A_n}_{G_n,A_n(l_1)}$ are
geometrically well-behaved. This was utilized in the proof of Lemma
3.2 to claim that the support of the barycentric extension of the
Whitehead homeomorphism is ``essentially'' at the place where the
Whitehead move exchanges diagonals.

\vskip .3 cm

\paragraph{\bf Remark} Note that $A_n\to id$ does not imply that $l$
and $A_n(l)$ are close uniformly for all edges $l$ of $\mathcal{F}$.
In fact, if the distance from $l$ to $l_0$ goes to infinity then the
distance between $l$ and $A_n(l)$ goes to infinity for each $n$
fixed. It is essential that we choose Whitehead moves along $G_n\{
l_1\}$ and $G_n\{ A_n(l_1)\}$ with $l_1$ close to $l_0$ and fixed.
Since $l_1$ and $A_n(l_1)$ are close, then their images under $G_n$
are close which allows us to compare the two maps along the orbits
$G_n\{ l_1\}$ and $G_n\{ A_n(l_1)\}$ whose corresponding elements
are close. The crucial fact that allows our method to work is that
maps $E(f_{id})$ and $E(g_{A_n})$ have small Beltrami dilatations
far away from the place where the Whitehead moves exchange the
diagonals because we do not have a uniform geometric control over
the maps away from the places where the diagonals are exchanged, see
above remark.

\vskip .3 cm

Let $K_N=\max\{ K_N(z_0,\mathcal{F}),K_N(z_0,\mathcal{F}^A)\}$.
Then, by (\ref{difference}) and by the choice of the above
neighborhoods of the $G_n$-orbit of $z_0$, $w_n$ belongs to the
$K_N$-neighborhood of the orbit $G_n\{ z_0\}$. Thus there exists
$\gamma_n\in G_n$ such that $w_n$ is in the $K_N$-neighborhood of
$\gamma_n(z_0)$ for each $n\in\mathbb{N}$. By the transitivity of
$PSL_2(\mathbb{Z})$ on the oriented edges of $\mathcal{F}$ and by
the transitivity of $A_n\circ PSL_2(\mathbb{Z})\circ A_n^{-1}$ on
the oriented edges of $\mathcal{F}^{A_n}$, there exist $\delta_n\in
PSL_2(\mathbb{Z})$ and $\delta_n'\in A_n\circ PSL_2(\mathbb{Z})\circ
A_n^{-1}$ such that $\delta_n\circ f_{id}\circ\gamma_n(l_0)=l_0$ and
$\delta_n'\circ g_{A_n}\circ\gamma_n(l_0)=l_0$.

\vskip .3 cm Since $A_n\to id$ as $n\to\infty$, it follows that
$[G_0:G_{A_n}]\to\infty$ as $n\to\infty$. This implies that
$[G_0:G_n]\to\infty$ as $n\to\infty$. The homeomorphism
$\delta_n\circ f_{id}\circ\gamma_n$ maps $\mathcal{F}_{G_n,l_1}$
onto $\mathcal{F}$, and the homeomorphism $\delta_n'\circ
g_{A_n}\circ\gamma_n$ maps $\mathcal{F}_{G_n,A_n(l_1)}^{A_n}$ onto
$\mathcal{F}^{A_n}$. The sequence of tessellations
$\mathcal{F}_{G_n,l_1}$ converges to the tessellation
$\mathcal{F}_{l_1}$ which differs from the Farey tessellation
$\mathcal{F}$ by the Whitehead move on the single edge $l_1$; the
sequence of tessellations $\mathcal{F}^{A_n}_{G_n,A_n(l_1)}$
converges to the tessellation $\mathcal{F}_{l_1}$ as well; and
$\mathcal{F}^{A_n}$ converges to the Farey tessellation
$\mathcal{F}$ (the convergence is in the Hausdorff topology on
compact subsets of the space of geodesics in $\D$). The above
convergence of the tessellations and the normalizations of
$\delta_n\circ f_{id}\circ\gamma_n$ and of $\delta_n'\circ
g_{A_n}\circ\gamma_n$ implies that both maps pointwise converge to
the Whitehead homeomorphism $f_{l_1}$ which maps $\mathcal{F}_{l_1}$
onto the Farey tessellation  $\mathcal{F}$ and which fixes $l_0$
(see Figure 4). \vskip .3 cm Since $\delta_n\circ
f_{id}\circ\gamma_n\to f_{l_1}$ and $\delta_n'\circ
g_{A_n}\circ\gamma_n\to f_{l_1}$ pointwise as $n\to\infty$, it
follows that the Beltrami dilatations $\mu (\delta_n\circ
E(f_{id})\circ\gamma_n)=\mu (E(f_{id})\circ\gamma_n)$ and $\mu
(\delta_n'\circ E(g_{A_n})\circ\gamma_n)=\mu
(E(g_{A_n})\circ\gamma_n)$ converge uniformly on compact subsets of
$\D$ to the Beltrami dilatation $\mu (E(f_{l_1}))$ of the
barycentric extension $E(f_{l_1})$ of $f_{l_1}$ (see \cite{DE}).
Since $\gamma_n^{-1}(w_n)$ belongs to the $K_N$-neighborhood of
$z_0$, it follows that
$$
|\mu (E(f_{id})\circ\gamma_n)(\gamma_n^{-1}(w_n))- \mu
(E(g_{A_n})\circ\gamma_n)(\gamma_n^{-1}(w_n))|\to 0
$$
as $n\to\infty$. This is the same as
$$
|\mu (E(f_{id}))(w_n)- \mu (E(g_{A_n}))(w_n)|\to 0
$$
as $n\to\infty$. But this is in the contradiction with
(\ref{difference}). The contradiction proves the theorem. $\Box$

\vskip .3 cm To finish the proof of  Theorem 1 it is remains to
establish that the circle homeomorphisms $g_A\circ f_{id}^{-1}$ are
not conformal maps. In fact, the proof below shows this for $A$
close enough to the identity and when the corresponding group $G$
(that determines the map $f_A$) has a sufficiently large index.

\vskip .3 cm

\paragraph{\bf Theorem 1}{\it For every $\epsilon >0$ there exist two
finite index subgroups of $PSL_2(\mathbb{Z})$ which are conjugated
by a $(1+\epsilon )$-quasisymmetric homeomorphism of the unit circle
and this conjugation homeomorphism is not conformal.}

\vskip .3 cm

\paragraph{\bf Proof} Recall that the Whitehead
homeomorphism $f_{id}$ maps $\mathcal{F}_{G,l_1}$ onto the Farey
tessellation $\mathcal{F}$, and that the Whitehead homeomorphism
$g_A$ maps $\mathcal{F}^A_{G,A(l_1)}$ onto $\mathcal{F}^A$, where
$G<G_A$ is any subgroup of finite index. By Theorem 3.3, there
exists a neighborhood $\mathcal{U}_{id}$ of the identity in
$PSL_2(\mathbb{Q})$ such that for any hyperbolic translation
$A\in\mathcal{U}_{id}$ whose oriented axis is $l_0$ the composition
$E(g_A)\circ E(f_{id})^{-1}$ of the barycentric extension of $g_{A}$
and $f_{id}$ has the quasiconformal constant less than $1+\epsilon$.
It is enough to show that $g_A\circ f_{id}^{-1}$ conjugates a finite
index subgroup of $G_A$ onto another finite index subgroup of
$PSL_2(\mathbb{Z})$ and that $g_A\circ f_{id}^{-1}$ is not conformal
(that is, the homeomorphism $g_A\circ f_{id}^{-1}$ is not a M\"obius
transformation).

\vskip .3 cm

Let $f_n$ and $g_n$ be two Whitehead homeomorphisms corresponding to
the Whitehead moves along the orbits $G_n\{ l_1\}$ and $G_n\{
A(l_1)\}$ on $\mathcal{F}$ and $\mathcal{F}^A$, where $G_n<G_A$ is a
sequence of finite index subgroups with $\cap_{n=1}^{\infty}G_n=\{
id\}$. In this case, the sequence $f_n$ converges pointwise to the
Whitehead homeomorphism $f_{l_1}$ which maps the tessellation
$\mathcal{F}_{l_1}$ onto the Farey tessellation $\mathcal{F}$, where
$\mathcal{F}_{l_1}$ is the image of the Farey tessellation
$\mathcal{F}$ under the Whitehead move on a single edge $l_1$. The
sequence $g_n$ pointwise converges to the Whitehead homeomorphism
$g_{A(l_1)}$ which maps $\mathcal{F}^A_{A(l_1)}$ onto
$\mathcal{F}^A$, where $\mathcal{F}^A_{A(l_1)}$ is the image of
$\mathcal{F}^A$ under the Whitehead move on a single edge $A(l_1)$
(see Figure 4).

\begin{figure}
\centering
\includegraphics[scale=0.4]{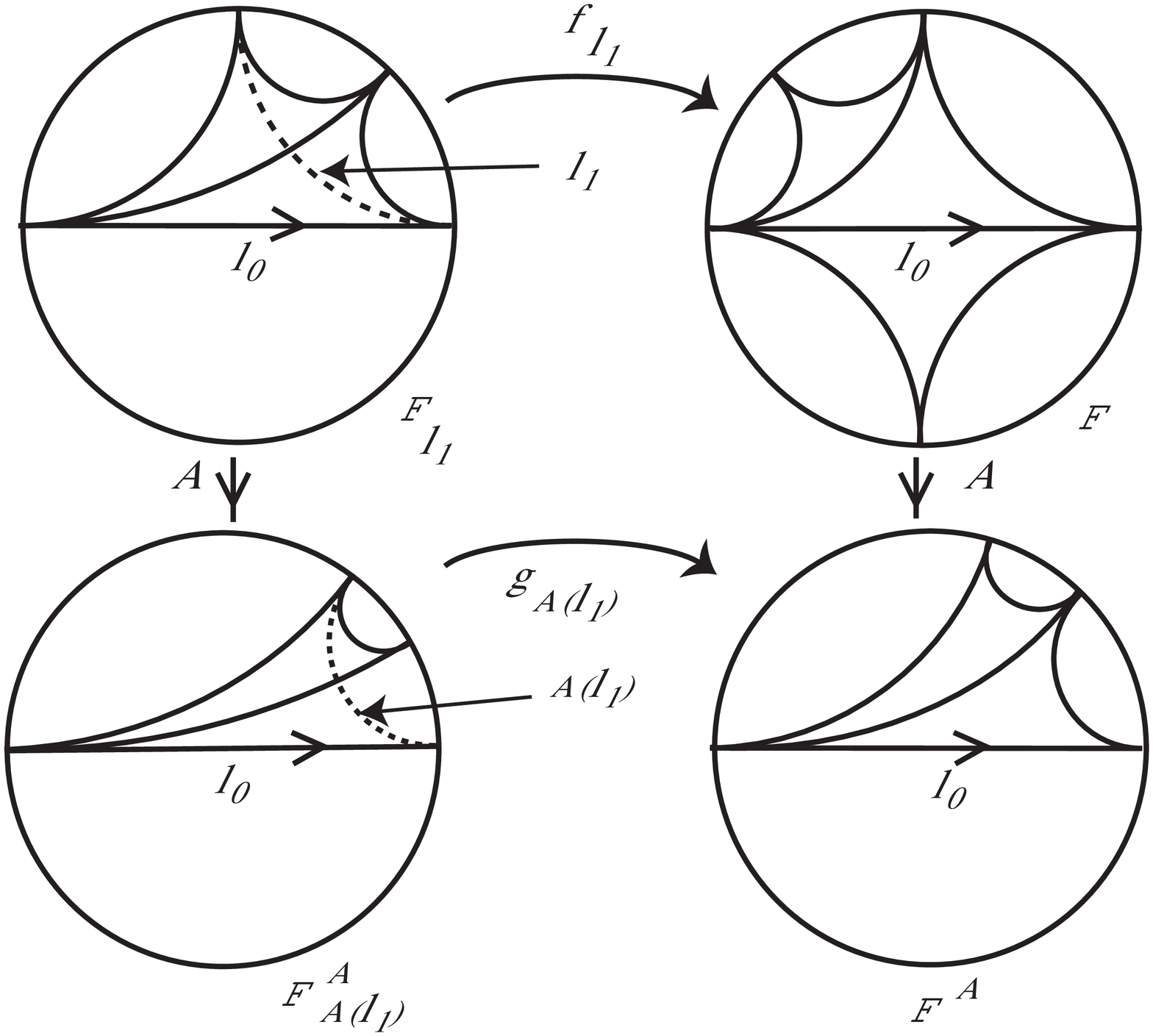}
\caption{} \label{}
\end{figure}

\vskip .3 cm

To see that $g_n\circ (f_n)^{-1}$ is not a M\"obius map for $n$
large enough, it is enough to show that $g_{A(l_1)}\circ
f_{l_1}^{-1}$ is not a M\"obius map. We note that the Whitehead
homeomorphisms $f_{l_1}$ is given by $f_{l_1}^{-1}=id$ on
$[-1,1]\subset S^1$, where

$$[-1,1]=\{ z\in S^1;\text{ $-1,z,1$ are
in the counterclockwise order}\}.
$$
\noindent The restriction $f_{l_1}^{-1}|_{[x_0,-1]}$ is the unique
element of $PSL_2(\mathbb{Z})$ which maps the oriented geodesic
$(-1,x_0)$ onto the oriented geodesic $(-1,i)$, where $x_0$ is the
third vertex of the complementary triangle of $\mathcal{F}$ to the
left of the oriented geodesic $(-1,i)$ with $(-1,i)$ on its
boundary. Also, the restriction $f_{l_1}^{-1}|_{[i,x_0]}$ is the
unique element of $PSL_2(\mathbb{Z})$ which maps the oriented
geodesic $(x_0,i)$ onto the oriented geodesic $(i,y_0)$, where $y_0$
is the third vertex of the complementary triangle of $\mathcal{F}$
to the left of $(i,1)$ with $(i,1)$ on its boundary. Finally,
$f_{l_1}^{-1}|_{[1,i]}$ is the unique element of $PSL_2(\mathbb{Z})$
which maps $(i,1)$ onto $(y_0,1)$ (see Figure 5). Thus, the
homeomorphism $f_{l_1}^{-1}$ is a piecewise $PSL_2(\mathbb{Z})$ with
four singular points $-1$, $1$, $i$ and $x_0$.  At these points the
map $f_{l_1}^{-1}$ changes its definition from one to another
element of $PSL_2(\mathbb{Z})$ (see Figure 5). It is interesting to
note (although we do not use this fact) that the homeomorphism
$f_{l_1}^{-1}$ is differentiable at every point on the circle.

\vskip .3 cm

Similarly, the singular points where $g_{A(l_1)}$ changes it
definition from one to another $PSL_2(\mathbb{Z})$ element are $-1$,
$1$, $A(y_0)$ and $A(i)$. Then $g_{A(l_1)}\circ f_{l_1}^{-1}$ is the
identity on $[-1,1]$, but at the point $i$ we have that
$f_{l_1}^{-1}$ changes its definition from one to another element of
$PSL_2(\mathbb{Z})$, while the restriction of $g_{A(l_1)}$ to a
neighborhood of $f_{l_1}^{-1}(i)=y_0$ equals a single element of
$PSL_2(\mathbb{Z})$ (because $A(y_0)\neq y_0$). This implies that
$g_{A(l_1)}\circ f_{l_1}^{-1}$ is not the identity in a neighborhood
of $i$. Thus $g_{A(l_1)}\circ f_{l_1}^{-1}$ is not a M\"obius map on
$S^1$. Consequently, $g_n\circ (f_n)^{-1}$ is not a M\"obius map for
all $n$ large enough. This completes the proof of Theorem 1. $\Box$

\begin{figure}
\centering
\includegraphics[scale=0.4]{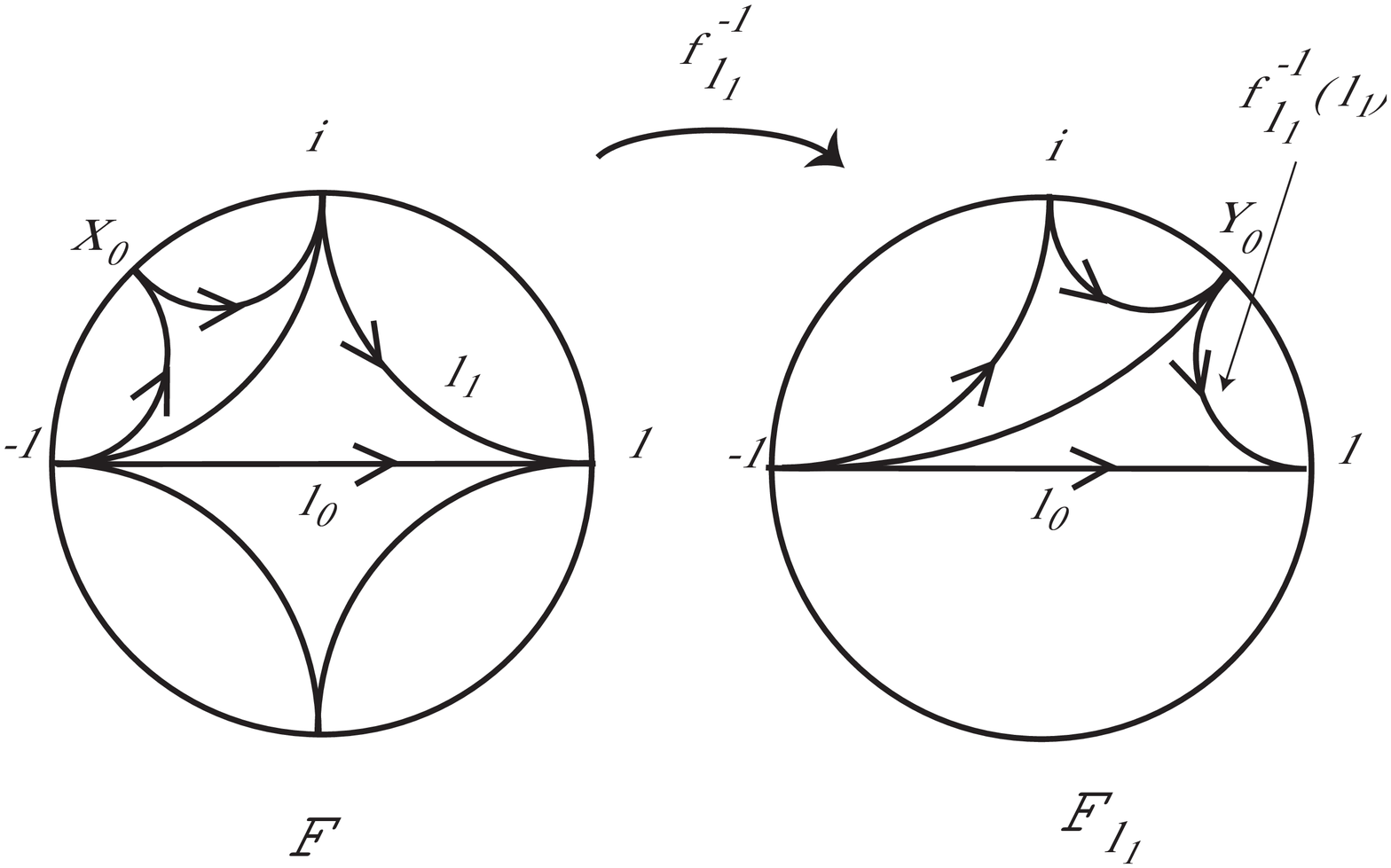}
\caption{} \label{}
\end{figure}

\vskip .3 cm

\section{The punctured solenoid $\S$}

Ehrenpreis conjecture asks whether any two compact Riemann surfaces
have finite regular covers which are close to being conformal, i.e.
if there exists a quasiconformal map between the covers which has
quasiconformal constant arbitrary close to $1$. Instead of taking
two arbitrary compact Riemann surfaces at a time and studying their
covers, an idea of Sullivan is to take all compact Riemann surfaces
at one time (i.e. in a single space) and keep track of the lifts via
the action of a Modular group. The same idea can be used for
punctured surfaces. We give more details below.

\vskip .3 cm

Let $T_0$ be the (once-punctured) Modular torus and let
$G_0<PSL_2(\mathbb{Z})$ be its universal covering group, i.e.
$T_0\equiv\D /G_0$. Let $S\to T_0$ be any finite regular covering of
$T_0$. Then there exists a natural isometric embedding
$T(T_0)\hookrightarrow T(S)$ of the Teichm\"uller space $T(T_0)$ of
the Modular torus $T_0$ into the Teichm\"uller space $T(S)$ of the
covering surface $S$. Moreover, if for a finite regular covering
$S_1\to T_0$ there exist finite regular coverings $S_2\to T_0$ and
$S_1\to S_2$ such that the composition $S_1\to S_2\to T_0$ is equal
to the original covering $S_1\to T_0$ then there is a natural
embedding $T(S_2)\hookrightarrow T(S_1)$ such that the image of
$T(T_0)$ in $T(S_2)$ is mapped onto the image of $T(T_0)$ in
$T(S_1)$. The inverse system of the finite regular coverings of
$T_0$ induces a direct system of Teichm\"uller spaces of the
covering surfaces. We denote the direct limit of the system of
Teichm\"uller spaces of all finite regular coverings of $T_0$ by
$T_{\infty}$ (see \cite{NS} for more details). The peripheral
preserving commensurator group $Comm_{per}(G_0)$ of the Modular
torus group $\pi_1(T_0)=G_0$ keeps track of different lifts of the
complex structure on the Modular torus. Thus, the Ehrenpreis
conjecture is equivalent to the statement whether $Comm_{per}(G_0)$
has dense orbit in $T_{\infty}$.

\vskip .3cm

The Teichm\"uller space $T(S)$ of any finite regular covering $S\to
T_0$ embeds in the universal Teichm\"uller space $T(\D)$ (i.e. the
Teichm\"uller space of the unit disk $\D$) as follows. Let $G<G_0$
be such that $\D /G\equiv S$ and that the covering $\D /G\to\D /G_0$
is conformally equivalent to $S\to T_0$. Then the image of $T(S)$ in
$T(\D )$ consists, up to an equivalence, of all Beltrami dilatations
$\mu$ on $\D$ such that
\begin{equation}
\label{inv}
 \mu (A(z))\frac{\overline{A'(z)}}{A'(z)}=\mu (z)
\end{equation}
for all $A\in G$ and $z\in\D$. Two Beltrami dilatations $\mu$ and
$\nu$ are equivalent if there is a quasiconformal map of $\D$ whose
Beltrami dilatation is $\mu -\nu$ and which extends to the identity
on ${\bf S}^1$.

\vskip .3 cm

Thus the image of the embedding $T_{\infty}\hookrightarrow T(\D )$
consists of all Beltrami dilatations $\mu$ on $\D$ which satisfy
(\ref{inv}) for some finite index subgroup $G$ of $G_0$. The image
of $T(S)$ under the embedding $T(S)\hookrightarrow T(\D )$ is a
finite-dimensional complex submanifold of $T(\D )$ but the embedding
is not an isometry for the Teichm\"uller metric (in fact, it is a
bi-biLipschitz map with the constant $1/3$ \cite{Mc}). The image of
$T_{\infty}$ in $T(\D )$ is not a closed subspace. The completion
$\overline{T_{\infty}}$ of the image of $T_{\infty}$ is a separable,
complex Banach submanifold of $T(\D )$ \cite{NS}. The completion
$\overline{T_{\infty}}$ consists of all Beltrami coefficients $\mu$
on $\D$ which are {\it almost invariant} under $G_0$ (modulo the
equivalence relation), i.e. $\overline{T_{\infty}}$ consists of all
$\mu$ which satisfy
$$
\sup_{A\in G_n}\|\mu\circ
A\frac{\overline{A'}}{A'}-\mu\|_{\infty}\to 0
$$
as $n\to\infty$ where $G_n$ is the intersection of all subgroups of
$G_0$ of index at most $n$. (Note that each $G_n$ is a finite index
subgroup of $G_0$ and that $\cap_{n=1}^{\infty}G_n=\{ id\}$.) The
Ehrenpreis conjecture is also equivalent to the question whether
$Comm_{per}(\S )$ has dense orbits in $\overline{T_{\infty}}$.

\vskip .3 cm

The points in $\overline{T_{\infty}}\setminus T_{\infty}$ are
obtained as limits of quasiconformal maps between finite Riemann
surfaces. These points are represented by Beltrami coefficients on
$\D$ with the additional property of being almost invariant.
Sullivan \cite{Sul} introduced a new object, called the universal
hyperbolic solenoid, on which these limit points appear in a
geometrically natural fashion as quasiconformal maps between the
universal hyperbolic solenoids. (Note that the quasiconformal maps
between finite surfaces lift to quasiconformal maps between the
universal hyperbolic solenoids as well.) We study the punctured
solenoid $\S$ which is the counter part of the universal hyperbolic
solenoid in the presence of punctures. We give the details below. An
important feature is that the Teichm\"uller space $T(\S )$ of the
punctured solenoid is naturally isometrically and bi-holomorphically
equivalent to $\overline{T_{\infty}}$.

\vskip .3 cm

We recall the definition and basic properties of the punctured
solenoid $\S$ \cite{PS}, which is an analogue in the presence of the
punctures of the universal hyperbolic solenoid introduced by
Sullivan \cite{Sul}. We keep the notation $T_0$ for the Modular once
punctured torus. Then $T_0$ is conformally identified with $\D/G_0$,
where $\D$ is the unit disk and $G_0<PSL_2(\mathbb{Z})$ is the
unique uniformizing subgroup. Consider the family of all finite
degree, regular coverings of $T_0\equiv\D/G_0$. The family is
inverse directed and the inverse limit $\S$ is called the {\it
punctured solenoid} (see \cite{PS}). The punctured solenoid $\S$ is
a non-compact space which is locally homeomorphic to a $2$-disk
times a Cantor set($\equiv$the transverse set); each path component,
called a {\it leaf}, is a simple connected $2$-manifold which is
dense in $\S$. $\S$ has one topological end which is homeomorphic to
the product of a horoball and the transverse set of $\S$ modulo
continuous action by a countable group (see \cite{PS}). A fixed leaf
of $\S$ is called the {\it baseleaf}. The punctured solenoid $\S$
has a natural projection $\Pi :\S\to T_0$ such that the restriction
to each leaf is the universal covering. The hyperbolic metric on
$T_0$ lifts to a hyperbolic metric on each leaf of $\S$ and the
lifted leafwise hyperbolic metric on $\S$ is locally constant in the
transverse direction. The punctured solenoid has a unique holonomy
invariant transverse measure (see \cite{Odd}). When the transverse
measure is coupled with the leafwise measure given by the hyperbolic
area on leaves, the resulting product measure is finite on $\S$.

\vskip .3 cm

We define an arbitrary {\it marked hyperbolic punctured solenoid}
$X$ to be a topological space locally homeomorphic to a $2$-disk
times a Cantor set with transversely continuous leafwise hyperbolic
metrics together with a homeomorphism $f:\S\to X$ which is
quasiconformal when restricted to each leaf and whose leafwise
Beltrami coefficients are continuous in the essential supremum norm
over the global leaves for the transverse variation (for more
details see \cite{PS}). A hyperbolic metric on any finite sheeted,
unbranched cover of $T_0$ gives a marked hyperbolic punctured
solenoid whose hyperbolic metric is transversely locally constant
for a choice of local charts, and any transversely locally constant
punctured solenoid arises as a lift of a hyperbolic metric on a
finite area punctured surface. We define the {\it Teichm\"uller
space} $T(\S )$ of the punctured solenoid $\S$ to be the space of
all marked hyperbolic punctured solenoids modulo an equivalence
relation. Two marked hyperbolic punctured solenoids $f_1:\S\to X_1$
and $f_2:\S\to X_2$ are equivalent if there exist an isometry
$c:X_1\to X_2$ such that the map $f_2^{-1}\circ c\circ f_1:\S\to\S$
is isotopic to the identity; the equivalence class of $f_1:\S\to X$
is denoted by $[f_1]$. The set of all marked transversely locally
constant hyperbolic punctured solenoids is dense in $T(\S )$ (see
\cite{Sul}, \cite{PS}). The {\it basepoint} of $T(\S )$ is the
equivalence class $[id:\S\to\S ]$ of the identity map.

\vskip .3 cm

The {\it modular group} $Mod(\S )$ (also called the {\it baseleaf
preserving mapping class group} $MCG_{BLP}(\S )$ in the literature
\cite{Odd}, \cite{PS}) of the punctured solenoid $\S$ consists of
homotopy classes of quasiconformal self-maps of $\S$ which preserve
the baseleaf. The restriction to the baseleaf of $Mod(\S )$ gives an
injective representation of $Mod(\S )$ into the group of the
quasisymmetric maps of $S^1$ (see \cite{Odd}). From now on, we
identify $Mod(\S )$ with this representation without further
mentioning. Then $Mod(\S )$ consists of all quasisymmetric maps of
$S^1$ which conjugate a finite index subgroup of $G_0$ onto (a
possibly different) finite index subgroup of $G_0$ such that
parabolic (peripheral) elements are conjugated onto parabolic
(peripheral) elements (see \cite{Odd}, \cite{PS}). In other words,
$Mod(\S )$ is isomorphic to the subgroup $Comm_{per}(G_0 )$ of the
abstract commensurator of $G_0$ consisting of all elements which
preserve parabolics. In particular, $Mod(\S )$ contains
$PSL_2(\mathbb{Q})$ and all lifts to the unit disk $\D$ of the
mapping class groups of the surfaces $\D/K$, where $K<G_0$ ranges
over all finite index subgroups. Recall that the Teichm\"uller space
$T(\S )$ embeds into the universal Teichm\"uller space $T(\D )$ by
restricting the leafwise quasiconformal homeomorphisms of $\S$ onto
variable solenoids to the baseleaf. From now on, we identify $T(\S
)$ with its image in $T(\D )$ under this embedding. Then the
Ehrenpreis conjecture is equivalent to the question whether $Mod(\S
)$ has dense orbits in $T(\S )$.

\vskip .3 cm

If the Ehrenpreis conjecture is correct then we show that for any
$\epsilon >0$ and for any finite Riemann surface there exist two
finite degree, regular covers and a $(1+\epsilon )$-quasiconformal
map between the covers which is not homotopic to a conformal map. We
remark that Theorem 1 establishes the existence of such covers for
the Modular punctured torus $T_0$ and any of its finite regular
covers (without the assumption that the Ehrenpreis conjecture is
correct) but it seems a difficult question to establish the
existence of such covers for an arbitrary punctured surface.

\vskip .3 cm

\paragraph{\bf Lemma 4.1} {\it Assume that the Ehrenpreis conjecture
is correct. Then for any $\epsilon >0$ and for any finite Riemann
surface there exist two finite degree, regular covers and a
$(1+\epsilon )$-quasiconformal map between the covers which is not
homotopic to a conformal map.}

\vskip .3 cm

\paragraph{\bf Proof} Since we assumed that the Ehrenpreis
conjecture is correct, we get that the orbits of $Mod(\S )$ are
dense in $T(\S )$. Let $S$ be an arbitrary finite area punctured
hyperbolic surface and let $f:S_0\to S$ be a quasiconfomal map from
a finite, unbranched covering surface $S_0$ of the Modular punctured
torus $T_0$ to the surface $S$. We note that the map $f:S_0\to S$
lifts to a map $\tilde{f}:\S\to X$ and that the equivalence class
$[\tilde{f}]$ is an element of the Teichm\"uller space $T(\S )$.
Then the orbit under $Mod(\S )$ of $[\tilde{f}]\in T(\S )$ is dense
and, in particular, it accumulates onto $[\tilde{f}]$. Let $g_n\in
Mod(\S )$ be a sequence such that $[\tilde{f}\circ g_n^{-1}]\to
[\tilde{f}]$ as $n\to\infty$. This implies that the Beltrami
dilatation of $\tilde{f}\circ g_n^{-1}\circ \tilde{f}^{-1}$ is
converging to zero as $n\to\infty$ but the dilatation of any
quasiconformal map homotopic to $\tilde{f}\circ g_n^{-1}\circ
\tilde{f}^{-1}$ is not equal to zero. The map $\tilde{f}\circ
g_n^{-1}\circ \tilde{f}^{-1}$ conjugates a finite index subgroup of
$\pi_1(S)$ to a different subgroup of $\pi_1(S)$. The two subgroups
of $\pi_1(S)$ are conjugated by a quasiconformal map with a small
Beltrami dilatation and therefore they establish the lemma. $\Box$

\vskip .3 cm

In our previous work \cite{MS}, we find an infinite family of orbits
with accumulation points outside the orbits. In particular, $T(\S
)/Mod(\S )$ is not a Hausdorff space. The points of the orbits are
non-transversely locally constant points in $T(\S )$ (i.e. they
correspond to points in $\overline{T_{\infty}}\setminus T_{\infty}$)
and elements of $Mod(\S )$ which give accumulation points are in
$G_0$. In this paper, we find accumulation points outside the orbit
of a transversely locally constant point in $T(\S )$ (i.e. a point
in $T_{\infty}$) corresponding to the basepoint $[id:\S\to\S ]$ for
the hyperbolic metric on $\S$ obtained by the lift of the hyperbolic
metric on $T_0\equiv \D/G_0$.

\vskip .3 cm

We show that the closure of the orbit under the modular group of the
basepoint $[id]\in T(\S )$ is strictly larger than the orbit and
that the closure is uncountable. We use Baire category theorem and
Theorem 3.3 together with fact that elements in Theorem 3.3 are not
M\"obius which is established in the course of the proof of Theorem
1.

\vskip .3 cm

\paragraph{\bf Corollary 2} {\it The closure in the Teichm\"uller metric
of the orbit under the modular group $Mod(\S )$ of the basepoint in
$T(\S )$ is strictly larger than the orbit. Moreover, the closure of
the orbit is an uncountable set without isolated points.}

\vskip .3 cm

\paragraph{\bf Remark} We showed in \cite{MS} that there is a set of
points in $T(\S )$ such that the closures of their orbits under the
modular group $Mod(\S )$ are strictly larger than the orbits. These
points were all non-transversely locally constant points in $T(\S
)$. The above corollary establishes that the orbit of the basepoint,
which is a transversely locally constant point, under the modular
group $Mod(\S )$ contains points outside the orbit. However, it is
still unknown whether any of the accumulation points of the orbit of
the basepoint is a transversely locally constant point in $T(\S )$.
This is equivalent to the question whether we can find an example of
two non-commensurate surfaces for which the Ehrenpreis conjecture is
correct.

\vskip .3 cm

\paragraph{\bf Proof} We use Baire category theorem. Assume on the
contrary that the closure of the orbit under $Mod(\S )$ of the
basepoint in $T(\S )$ is equal to the orbit.

\vskip .3 cm

Thus the orbit is a closed subset in $T(\S )$, hence it is of the
second kind in itself (in the sense of Baire). We claim that there
exists a point of the orbit which is an isolated point. If not, then
each point of the orbit is nowhere dense. Since a single point in a
metric space is always a closed subset, it follows that the orbit
can be written as a countable union of its singletons (which are
nowhere dense closed sets). This contradicts the Baire theorem.

\vskip .3 cm

Therefore, at least one point $[f]\in T(\S )$ where $f\in Mod(\S )$
is isolated. Choose a sequence $f_n\in Mod(\S )$ satisfying the
properties in Theorem 3.3 such that $\|\mu (E(f_n))\|_{\infty}\to 0$
as $n\to\infty$. Then $f_n\circ f\in Mod(\S )$ is in the orbit of
the basepoint and $[f_n\circ f]\to [f]$ as $n\to\infty$ in the
Teichm\"uller metric on $T(\S )$. This is a contradiction.
Therefore, the closure of the orbit under $Mod(\S )$ of the
basepoint in $T(\S )$ is strictly larger than the orbit.

\vskip .3 cm

We proceed to prove that the closure of the orbit is uncountable.
Assume on the contrary that the orbit is countable. Then there
exists an isolated point $f$ of the closure of the orbit by the
above argument. The isolated point $f$ is necessarily in $Mod(\S )$
because the accumulation points of $Mod(\S )$ in $T(\S )$ are not
isolated. Then the above argument establishes a contradiction.
Therefore the closure is uncountable. $\Box$

\vskip .3 cm

In the corollary below, we show how to explicitly construct
sequences in the orbit under $Mod(\S )$ of the basepoint in $T(\S )$
which accumulate to points outside the orbit.

\vskip .3 cm

\paragraph{\bf Corollary 4.2} {\it There exists a sequence $f_n\in
Mod(\S )$ whose elements are constructed as in Theorem 3.3 such that
$[f_n\circ f_{n-1}\circ\cdots\circ f_1]$ converges to $[f]\in T(\S
)$, where $f\notin Mod(\S )$.}

\vskip .3 cm

\paragraph{\bf Proof} We choose a
sequence $f_n\in Mod(\S )$ such that $\|\mu (E(f_n))\|_{\infty}<
1/2^{n+1}$ and $\inf_{\{ g:g|_{S^1}= f_n|_{S^1}\}}\|\mu
(g)\|_{\infty}>0$ which is possible by Theorem 3.3, where $E(f_n)$
is the barycentric extension of $f_n$. Then the sequence
$E(f_n)\circ E(f_{n-1})\circ\cdots\circ E(f_1)$ has uniformly
convergent Beltrami coefficients. Therefore, $f_n\circ
f_{n-1}\circ\cdots\circ f_1$ converges in the Teichm\"uller metric
in $T(\S )$.

\vskip .3 cm

The sequence $f_n$ can be chosen such that each $f_n$ conjugates a
maximal finite index subgroup $H_n$ of $PSL_2(\mathbb{Z})$ onto
another finite index subgroup $K_n$ of $PSL_2(\mathbb{Z})$, where
$[PSL_2(\mathbb{Z}):H_n]\to\infty$ as $n\to\infty$. To see this we
take $A_n\to id$ as in Theorem 3.3 and for each $n$ consider a
sequence of finite index subgroups $G_{k_n}$ of $G_{A_n}$ which are
obtained as intersections of all subgroups of $G_{A_n}$ of index at
most $k_n$. We define $f_{n,k_n}$ using $G_{k_n}$ and $A_n$ as in
Theorem 3.3. Assume that for a fixed $n$ all such obtained maps
$f_{n,k_n}$ conjugate a fixed finite index subgroup
$K<PSL_2(\mathbb{Z})$ onto another (possibly different but of the
same index) subgroup of $PSL_2(\mathbb{Z})$. Since all $f_{n,k_n}$
fix $l_0$ by construction, it follows that $f_{n,k_n}$ converges
pointwise to $f_{A_n(l_1)}\circ f_{l_1}^{-1}$ as $k_n\to\infty$ and
$n$ fixed, where $f_{l_1}$ is a Whitehead homeomorphisms which maps
$\mathcal{F}_{l_1}$ onto $\mathcal{F}$, and $f_{A_n(l_1)}$ is a
Whitehead homeomorphism which maps $\mathcal{F}^{A_n}_{A_n(l_1)}$
onto $\mathcal{F}^{A_n(l_1)}$. We showed in the proof of Theorem 1
that $f_{A_n(l_1)}\circ f_{l_1}^{-1}$ is a piecewise M\"obius but
not a M\"obius map. Then $f_{A_n(l_1)}\circ f_{l_1}^{-1}$ cannot
conjugate a finite index subgroup of $PSL_2(\mathbb{Z})$ onto
itself. Therefore, for $k_n$ large enough, $f_{n,k_n}$ does not
conjugate $K$ onto another finite index subgroup. Therefore, we can
choose $k_n$ large enough such that $f_{n,k_n}$ conjugates a maximal
subgroup $H_n$ of finite index in $PSL_2(\mathbb{Z})$ onto another
finite index subgroup with $[PSL_2(\mathbb{Z}):H_n]\to\infty$ as
$n\to\infty$. Moreover, the subgroups $H_n$ can be chosen such that
$\cap_{n=1}^{\infty}H_n=\{ id\}$. If we further enlarge $k_n$ the
above remains true.

\vskip .3 cm

We put the elements of $G_0$ in a sequence
$\{\gamma_i\}_{i=1}^{\infty}$. For each $n$ and for each $i$, there
exists $k_n(i)$ such that $\gamma_i\notin G_{k_n(i)}$. This implies
that $\gamma_i\notin G_j$ for all $j\geq k_n(i)$. We apply Cantor's
diagonal argument. For $\gamma_1$, we choose $n=n(1)$ arbitrary.
Then $f_{n(1),k_n(1)}$ does not conjugate $\gamma_1$ onto any other
element of $PSL_2(\mathbb{Z})$. We consider the minimal dilatation
$$
\min (f_{n(1),k_n(1)}) (\gamma_1):=\inf_{\{
g;g|_{S^1}=f_{n(1),k_n(1)}\}}\|\mu (g)-\mu
(g\circ\gamma_1)\|_{\infty}.
$$

 In the case that $\min (f_{n(1),k_n(1)})(\gamma_2)=0$ then
we choose $n(2)$ small enough such that $\|\mu
(E(f_{n(2),k_n(2)}))\|_{\infty}<\frac{1}{4}\min (f_{n(1),k_n(1)})
(\gamma_1)$. If $\min (f_{n(1),k_n(1)})(\gamma_2)>0$ we choose
$n(2)$ such that $\|\mu (E(f_{n(2),k_n(2)}))\|_{\infty}$ is less
than $\frac{1}{4}$ of the minimum of $\min (f_{n(1),k_n(1)})
(\gamma_1)$ and $\min ( f_{n(1), k_n(1)}) (\gamma_2)$. In both cases
we are guaranteed that $f_{n(2),k_n(2)}\circ f_{n(1),k_n(1)}$ does
not conjugate $\gamma_1$ or $\gamma_2$ onto an element of
$PSL_2(\mathbb{Z})$. For $\gamma_3$, we choose $n(3)$ such that
$\|\mu (E(f_{n(3),k_n(3)}))\|_{\infty}$ is less than $\frac{1}{4}$
of the minimum of $\min (f_{n(2),k_n(2)}) (\gamma_1)$, $\min
(f_{n(2),k_n(2)}\circ f_{n(1),k_n(1)})(\gamma_2)$ and $\min
(f_{n(2),k_n(2)}\circ f_{n(1),k_n(1)})(\gamma_3)$. This guarantees
that $f_{n(3),k_n(3)}\circ f_{n(2),k_n(2)}\circ f_{n(1),k_n(1)}$
does not conjugate $\gamma_i$, for $i=1,2,3$ onto elements of
$PSL_2(\mathbb{Z})$. We continue this process for all
$i\in\mathbb{N}$.

\vskip .3 cm

By our choice of $n(i)$, the series $\sum_{i=1}^{\infty}\|\mu
(E(f_{n(i),k_n(i)}))\|_{\infty}$ converges. Thus the sequence
$f_{n(i),k_n(i)}\circ f_{n(i-1),k_n(i-1)}\circ\cdots \circ f_1$
converges in the Teichm\"uller metric. By the above choices, the
limit does not conjugate a single element of $G_0$ onto any other
element of $PSL_2(\mathbb{Z})$. Thus the limit is not in the orbit
under the modular group $Mod(\S )$ of the base point in $T(\S )$.
$\Box$

\end{document}